\documentclass[a4paper, 11pt, twoside]{article}
\usepackage{mfpaperstuff}
\usetikzlibrary{matrix,arrows,decorations.pathmorphing,decorations.pathreplacing}
\usetikzlibrary{decorations.markings}

\DeclareMathOperator{\Hom}{Hom}
\DeclareMathOperator{\codeg}{codeg}
\DeclareMathOperator{\sgn}{sgn}
\DeclareMathOperator{\DHom}{DHom}
\DeclareMathOperator{\dhom}{DHom}

\DeclareMathOperator{\res}{res}
\DeclareMathOperator{\cont}{cont}
\DeclareMathOperator{\lab}{lab}

\begin{document}
\newcommand\bil[2]{(#1\,|\,#2)}
\newcommand\cstr{_{\mathrm{c}}}
\renewcommand\hom{\operatorname{Hom}}
\newcommand\shom[2]{\Hom_{\hhh n}(\spe{#1},\spe{#2})}
\newcommand\rshom[2]{\Hom_{\hhh n}(\rspe{#1},\rspe{#2})}
\newcommand\chom[2]{\DHom_{\hhh n}(\spe{#1},\spe{#2})}
\newcommand\rhom[2]{\DHom_{\hhh n}(\rspe{#1},\rspe{#2})}
\renewcommand\unlhd\domby
\renewcommand\unrhd\dom
\renewcommand{\crefrangeconjunction}{--}
\medmuskip=3mu plus 1mu minus 1mu
\newcommand\mo{{-}1}
\newcommand\yfz{\Yfillopacity0}
\newcommand\yfo{\Yfillopacity1}
\newcommand\tpf{2.5}
\newcommand\ten{10}
\newcommand\eleven{11}
\newcommand\twelve{12}
\newcommand\thirteen{13}
\newcommand\fourteen{14}
\newcommand\hf{.5}
\newcommand\sq{1.5}
\newcommand\one{\mathbbm{1}}
\newcommand\bsm{\begin{smallmatrix}}
\newcommand\esm{\end{smallmatrix}}
\newcommand{\rt}[1]{\rotatebox{90}{$#1$}}
\newcommand\id{\operatorname{id}}
\newcommand{\ol}{\overline}
\newcommand\partn{\mathcal{P}}
\newcommand{\nchar}{\operatorname{char}}
\newcommand{\thmcite}[2]{\textup{\textbf{\cite[#2]{#1}}}\ }
\newcommand\zez{\mathbb{Z}/e\mathbb{Z}}
\newcommand{\hhh}{\mathscr{H}_}
\newcommand{\sect}[1]{\section{#1}}
\newcommand\cgs\succcurlyeq
\newcommand\cls\preccurlyeq
\newcommand\cg\succ
\newcommand\cl\prec
\newcommand\cnd{(\textasteriskcentered)}
\newcommand\cne{(\textdagger)}
\newcommand\trinom[3]{\big[#1\big|#2\big|#3\big]}
\newcommand\spe[1]{\operatorname{S}_{#1}}
\newcommand\rspe[1]{\operatorname{S}^{#1}}
\newcommand\spek[2]{\spe{#1|#2}}
\newcommand\rspek[2]{\rspe{#1|#2}}
\newcommand\std[1]{\operatorname{Std}(#1)}
\newcommand\cstd[2]{\operatorname{Std}_{#1}(#2)}
\newcommand\rstd[2]{\operatorname{Std}^{#1}(#2)}
\newcommand\stdlr[1]{\operatorname{Std}_{\mathrm{LR}}(#1)}
\newcommand\lf{_{\mathrm{L}}}
\newcommand\rg{_{\mathrm{R}}}
\newcommand\tp{_{\mathrm{T}}}
\newcommand\bt{_{\mathrm{B}}}
\newcommand\Lf{_{\mathrm{L}}}
\newcommand\cj{\#}
\newcommand\rj{\bar{\#}}
\newcommand\lr[2]{#1\lf\cj#2\rg}
\newcommand\ann{\operatorname{Ann}}
\newcommand\tabupto[2]{#1_{\downarrow#2}}
\newcommand\shp[2]{\operatorname{Shape}(\tabupto{#1}{#2})}
\newcommand\emc{$e$-multicharge\xspace}
\newcommand\emcs{$e$-multicharges\xspace}
\renewcommand\phi\varphi
\newcommand\dual{^\circledast}
\newcommand\df{\operatorname{def}}
\newcommand\mptn[2]{\scrp^{#1}_{#2}}
\newcommand\shift[1]{\operatorname{shift}_{#1}}
\newcommand\garn[1]{\mathsf g_{#1}}
\newcommand\rgarn[1]{\mathsf g^{#1}}
\Yvcentermath0
\Yboxdim{12pt}

\title{Generalised column removal for graded homomorphisms between Specht modules}
\author{Matthew Fayers \& Liron Speyer\\{\normalsize Queen Mary University of London, Mile End Road, London E1 4NS, U.K.}\\\texttt{\normalsize m.fayers@qmul.ac.uk}, \texttt{\normalsize l.speyer@qmul.ac.uk}}
\renewcommand\auth{Matthew Fayers \& Liron Speyer}
\runninghead{Generalised column removal for graded homomorphisms}

\toptitle

\begin{abstract}
Let $n$ be a positive integer, and let $\hhh n$ denote the affine KLR algebra in type A. Kleshchev, Mathas and Ram have given a homogeneous presentation for graded column Specht modules $\spe\la$ for $\hhh n$. Given two multipartitions $\la$ and $\mu$, we define the notion of a \emph{dominated} homomorphism $\spe\la\to\spe\mu$, and use the KMR presentation to prove a generalised column removal theorem for graded dominated homomorphisms between Specht modules. In the process, we prove some useful properties of $\hhh n$-homomorphisms between Specht modules which lead to an immediate corollary that, subject to a few demonstrably necessary conditions, every homomorphism $\spe\la\to\spe\mu$ is dominated, and in particular $\shom\la\mu=0$ unless $\la$ dominates $\mu$.

Brundan and Kleshchev show that certain cyclotomic quotients of $\hhh n$ are isomorphic to (degenerate) cyclotomic Hecke algebras of type A. Via this isomorphism, our results can be seen as a broad generalisation of the column removal results of Fayers and Lyle and of Lyle and Mathas; generalising both into arbitrary level and into the graded setting.
\end{abstract}

\section{Introduction}

The KLR algebras, or quiver Hecke algebras, were constructed independently by Khovanov and Lauda \cite{kl} and by Rouquier \cite{rouq}, and have since received an abundance of interest. This is, in some part, due to the powerful result of Brundan and Kleshchev in \cite{bk} that every (degenerate) Ariki--Koike algebra is isomorphic to a so-called cyclotomic quotient of a KLR algebra. The KLR algebras and their cyclotomic quotients are graded, and this allows us to study the graded representation theory of (degenerate) Ariki--Koike algebras, and in particular the graded representation theory of the symmetric groups. This motivates the study of KLR algebras, and in particular the study of their graded Specht modules. These were defined by Brundan, Kleshchev and Wang \cite{bkw}, and developed further by Kleshchev, Mathas and Ram \cite{kmr}, who gave a homogeneous presentation for each Specht module.

In trying to understand the (graded) structure of the Specht modules, the (graded) homomorphism spaces $\shom\la\mu$ are of particular interest. In the ungraded setting, these homomorphism spaces have received a great deal of attention in recent years. We concentrate in particular here on the row and column removal theorems for homomorphisms, proved by the first author and Lyle \cite{fl} for the symmetric group and generalised to Hecke algebras of type $A$ by Lyle and Mathas. In this paper, we provide graded versions of these theorems, while at the same time generalising them to higher levels so that they apply to all (degenerate) Ariki--Koike algebras.

In fact, our results apply not to all homomorphisms between two given Specht modules but only to those of a certain type, which we call \emph{dominated} homomorphisms. However, in many cases (for example, for the symmetric group in odd characteristic) every homomorphism between two Specht modules is dominated, so our results apply generally; in particular, via the Brundan--Kleshchev isomorphism mentioned above, we recover the original row and column removal theorems of Lyle and Mathas.

We now summarise the structure of the paper. In \cref{backsec}, we introduce the combinatorics necessary for our purposes, as well as the set-up of the KLR algebras and their Specht modules. We proceed in \cref{domsec} by introducing \emph{dominated tableaux} and the corresponding dominated homomorphisms. \cref{crhsec} gives our main results pertaining to generalised column removal for homomorphisms. Finally, in \cref{indsec} we provide an index of notation for the reader's convenient reference.

\begin{acks}
The second author would like to thank Queen Mary University of London, without whose funding this work would not have been possible. The second author must also thank Professor Andrew Mathas, at the University of Sydney, with whom this work began. The visit to the University of Sydney was funded by the Eileen Colyer Prize and the Australian Research Council grant DP110100050 ``Graded representations of Hecke algebras''. The authors thank the anonymous referee for his or her extensive and helpful comments.
\end{acks}

\section{Background}\label{backsec}

In this section we recall some background and set up some notation. This varies from \cite{kmr} in only a few details.

\subsection{The symmetric group}\label{symsec}

Let $\sss n$ denote the symmetric group of degree $n$. Let $s_1,\dots,s_{n-1}$ denote the standard Coxeter generators of $\sss n$, i.e.\ $s_i$ is the transposition $(i,i+1)$. Given $w\in\sss n$, a \emph{reduced expression} for $w$ is an expression $w=s_{i_1}\dots s_{i_l}$ with $l$ as small as possible; we call $l=l(w)$ the \emph{length} of $w$.

We will need to use two natural partial orders on $\sss n$. If $w,x\in\sss n$, then we say that $x$ is smaller than $w$ in the \emph{left order} (and write $x\ls\Lf w$) if $l(w)=l(wx^{-1})+l(x)$; this is equivalent to the statement that there is a reduced expression for $w$ which has a reduced expression for $x$ as a suffix.

More important will be the \emph{Bruhat order} on $\sss n$: if $w,x\in\sss n$, then we say that $x$ is smaller than $w$ in the Bruhat order (and write $x\cls w$) if there is a reduced expression for $w$ which has a (possibly non-reduced) expression for $x$ as a subsequence. In fact \cite[Theorem 5.10]{hum}, if $x\cls w$, then every reduced expression has a reduced expression for $x$ as a subsequence.

The following proposition gives an alternative characterisation of the Bruhat order.

\begin{propnc}{hum}{\S5.9}\label{tbruhat}
Suppose $w,x\in\sss n$. Then $w\cls x$ if and only if there are $w=w_0,w_1,\dots,w_r=x$ such that for each $1\ls i\ls r$ we have $w_i=(u_i, v_i)w_{i-1}$, where $1\ls u_i<v_i\ls n$ and $w_{i-1}^{-1}(u_i)<w_{i-1}^{-1}(v_i)$.
\end{propnc}

Later we shall need the following lemma; in fact, this is a special case of Deodhar's `property Z' \cite[Theorem 1.1]{deod}.

\begin{lemma}\label{brulem}
Suppose $w,x\in\sss n$ with $x\cl w$. If $l(s_iw)<l(w)$ while $l(s_ix)>l(x)$, then $s_ix\cls w$.
\end{lemma}

\begin{pf}
Since $l(s_iw)<l(w)$, $w$ has a reduced expression $s$ beginning with $s_i$. We can find a reduced expression for $x$ as a subexpression of $s$, and this subexpression cannot include the first term $s_i$, since $l(s_ix)>l(x)$. So we can add the initial $s_i$ to the subexpression to get a reduced expression for $s_ix$ as a subexpression of $s$.
\end{pf}

We end this subsection by defining some very natural and useful homomorphisms. Suppose $1\leq m \leq n$ and $0\leq k\leq n-m$, and define the homomorphism $\shift k:\sss m \rightarrow \sss n$ by $s_i\mapsto s_{i+k}$ for every $i$. Note that if $k=0$, this is the natural embedding.

\subsection{Lie-theoretic notation}\label{liesec}

Throughout this paper $e$ is a fixed element of the set $\{2,3,4,\dots\}\cup\{\infty\}$. If $e=\infty$ then we set $I:=\bbz$, while if $e<\infty$ then we set $I:=\zez$; we may identify $I$ with the set $\{0,\dots,e-1\}$ when convenient. The Cartan matrix $(a_{i,j})_{i,j\in I}$ is defined by $a_{ij}=2\delta_{ij}-\delta_{i(j+1)}-\delta_{i(j-1)}$.

Let $\Gamma$ be the quiver with vertex set $I$ and an arrow from $i$ to $i-1$ for each $i$. (Note that this convention is the same as that in \cite{kmr}, and opposite to that in \cite{bk,bkw}.) The quiver $\Gamma$ is pictured below for some values of $e$.
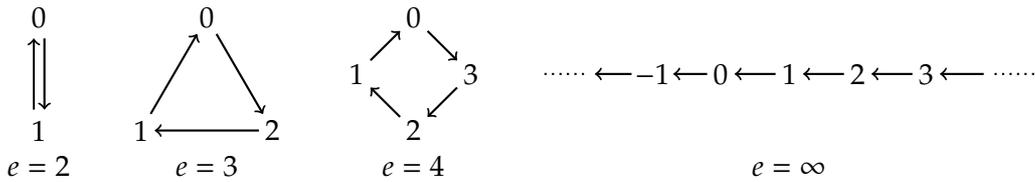
\[
\begin{array}{c@{\qquad}c@{\qquad}c@{\qquad}c}
\begin{tikzpicture}[thick,baseline=-.75cm]
\foreach\x in{0,1}
{\node(\x)at(0,-1.5*\x){$\mathclap{\x}$};
\node(r\x)at(.07,-1.5*\x){$\phantom{\x}$};
\node(l\x)at(-.07,-1.5*\x){$\phantom{\x}$};}
\draw[->](r0)--(r1);
\draw[->](l1)--(l0);
\end{tikzpicture}
&
\begin{tikzpicture}[thick,baseline=.25cm]
\foreach\x in{0,1,2}
{\node[inner sep=2pt](\x)at(90+\x*120:1){$\x$};
\node[inner sep=2pt](p\x)at(330+\x*120:1){$\phantom{\x}$};
\draw[->](\x)--(p\x);}
\end{tikzpicture}
&
\begin{tikzpicture}[thick,scale=.75,baseline=0pt]
\foreach\x in{0,1,...,3}
{\node[inner sep=2pt](\x)at(90+\x*90:1){$\x$};
\node[inner sep=2pt](p\x)at(\x*90:1){$\phantom{\x}$};
\draw[->](\x)--(p\x);}
\end{tikzpicture}
&
\begin{tikzpicture}[thick,scale=.9,baseline=0pt]
\foreach\x in{0,1,2,3}
\node[inner sep=2pt](\x)at(\x,0){$\x$};
\foreach\x in{-1}
\node[inner sep=1pt](\x)at(\x,0){$\x$};
\foreach\x in{-2,4}
\node(\x)at(\x,0){};
\draw[->](4)--(3);
\draw[->](3)--(2);
\draw[->](2)--(1);
\draw[->](1)--(0);
\draw[->](0)--(-1);
\draw[->](-1)--(-2);
\draw[dotted](4,0)--++(.6,0);
\draw[dotted](-2,0)--++(-.6,0);
\end{tikzpicture}
\\
e=2&e=3&e=4&e=\infty
\end{array}
\]
In the relations we give below, we use arrows with reference to $\Gamma$; thus we may write $i\to j$ to mean that $e\neq2$ and $j=i-1$, or $i\rightleftarrows j$ to mean that $e=2$ and $j=i-1$.

We adopt standard notation from Kac's book \cite{kac} for the Kac--Moody algebra associated to the Cartan matrix $(a_{i,j})_{i,j\in I}$; in particular, we have fundamental dominant weights $\La_i$ and simple roots $\alpha_i$ for $i\in I$, and an invariant symmetric bilinear form $\bil\,\,$ satisfying $\bil{\La_i}{\alpha_j}=\delta_{ij}$ and $\bil{\alpha_i}{\alpha_j}=a_{ij}$ for $i,j\in I$. A \emph{root} is a linear combination $\sum_{i\in I}c_i\alpha_i$ with $c_i\in\bbz$ for each $i$; the \emph{height} of this root is defined to be $\sum_{i\in I}c_i$. Given roots $\alpha=\sum_{i\in I}c_i\alpha_i$ and $\beta=\sum_{i\in I}d_i\alpha_i$, we write $\alpha\gs\beta$ if $c_i\gs d_i$ for each $i$; we say that $\alpha$ is \emph{positive} if $\alpha\gs0$.

Let $I^l$ denote the set of all $l$-tuples of elements of $I$. We call an element of $I^l$ an \emph{$e$-multicharge of level $l$}. The symmetric group $\sss l$ acts on $I^l$ on the left by place permutations. Given an $e$-multicharge $\kappa=(\kappa_1,\dots,\kappa_l)$, we define a corresponding dominant weight $\La_\kappa:=\La_{\kappa_1}+\dots+\La_{\kappa_l}$. We then define the \emph{defect} of a root $\alpha$ (with respect to $\kappa$) to be
\[
\df(\alpha)=\bil{\La_\kappa}\alpha-\tfrac12\bil\alpha\alpha.
\]

\subsection{Multicompositions and multipartitions}\label{mptnsec}

A \emph{composition} is a sequence $\la=(\la_1,\la_2,\dots)$ of non-negative integers such that $\la_i=0$ for sufficiently large $i$. We write $|\la|$ for the sum $\la_1+\la_2+\cdots$. When writing compositions, we may omit trailing zeroes and group equal parts together with a superscript. We write $\varnothing$ for the composition $(0,0,\dots)$. A \emph{partition} is a composition $\la$ for which $\la_1\gs\la_2\gs\cdots$.

Now suppose $l\in\bbn$. An \emph{$l$-multicomposition} is an $l$-tuple $\la=(\la^{(1)},\dots,\la^{(l)})$ of compositions, which we refer to as the \emph{components} of $\la$. We write $|\la|=|\la^{(1)}|+\dots+|\la^{(l)}|$, and say that $\la$ is an $l$-multicomposition of $|\la|$. If the components of $\la$ are all partitions, then we say that $\la$ is an \emph{$l$-multipartition}. We write $\mptn ln$ for the set of $l$-multipartitions of $n$. We abuse notation by using $\varnothing$ also for the multipartition $(\varnothing,\dots,\varnothing)$.

If $\la$ and $\mu$ are $l$-multicompositions of $n$, then we say that $\la$ \emph{dominates} $\mu$, and write $\la\dom\mu$, if
\[
|\la^{(1)}|+\dots+|\la^{(m-1)}|+\la^{(m)}_1+\dots+\la^{(m)}_r\gs|\mu^{(1)}|+\dots+|\mu^{(m-1)}|+\mu^{(m)}_1+\dots+\mu^{(m)}_r
\]
for all $1\ls m\ls l$ and $r\gs0$.

If $\la$ is an $l$-multicomposition, the \emph{Young diagram} $[\la]$ is defined to be the set
\[
\rset{(r,c,m)\in\bbn\times\bbn\times\{1,\dots,l\}}{c\ls\la^{(m)}_r}.
\]
We refer to the elements of $[\la]$ as the \emph{nodes} of $\la$. We may also refer to $(r,c,m)$ as the $(r,c)$-node of $\la^{(m)}$. If $\la\in\mptn l n$, a node of $\la$ is \emph{removable} if it can be removed from $[\la]$ to leave the Young diagram of a smaller $l$-multipartition, while a node not in $[\la]$ is \emph{addable} if it can be added to $[\la]$ to leave the Young diagram of an $l$-multipartition.

We adopt an unusual (but in our view, extremely helpful) convention for drawing Young diagrams. We draw the nodes of each component as boxes in the plane, using the English convention, where the first coordinate increases down the page and the second coordinate increases from left to right. Then we arrange the diagrams for the components \emph{in a diagonal line from top right to bottom left}. For example, if $\la=\left((2^2),(2,1^2),(3,2)\right)\in\mptn3{13}$, then $[\la]$ is drawn as follows.
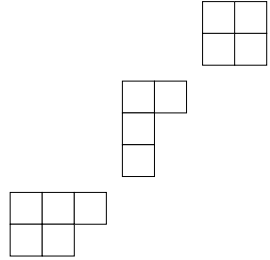
\[
\gyoung(^6;;,^6;;,/\hf,^3^\hf;;,^3^\hf;,^3^\hf;,/\hf,;;;,;;)
\]
We shall use directions such as left and right with reference to this convention; for example, we shall say that a node $(r,c,m)$ lies to the left of $(r',c',m')$ if either $m>m'$ or ($m=m'$ and $c<c'$). Similarly, we say that $(r,c,m)$ is above, or higher than, $(r',c',m')$ if either $m<m'$ or ($m=m'$ and $r<r'$).

If $\la$ is a partition, the \emph{conjugate partition} $\la'$ is defined by
\[
\la'_i=\left|\lset{j\gs1}{\la_j\gs i}\right|.
\]
If $\la$ is an $l$-multipartition, then the conjugate multipartition $\la'$ is given by
\[
\la'=({\la^{(l)}}',\dots,{\la^{(1)}}').
\]
Observe that with our convention, the Young diagram $[\la']$ may be obtained from $[\la]$ by reflecting in a diagonal line running from top left to bottom right.

\subsection{Tableaux}\label{tabsec}

If $\la\in\mptn ln$, a \emph{$\la$-tableau} is a bijection $\ttt:[\la]\to\{1,\dots,n\}$. We depict a $\la$-tableau $\ttt$ by drawing the Young diagram $[\la]$ and filling each box with its image under $\ttt$. $\ttt$ is \emph{row-strict} if its entries increase from left to right along each row of the diagram, and \emph{column-strict} if its entries increase down each column. $\ttt$ is \emph{standard} if it is both row- and column-strict. We write $\std\la$ for the set of standard $\la$-tableaux.

If $\ttt$ is a $\la$-tableau, then we define a $\la'$-tableau $\ttt'$ by
\[
\ttt'(r,c,m)=\ttt(c,r,l+1-m)
\]
for all $(r,c,m)\in[\la']$.

We import and modify some notation from \cite{bkw} and \cite{kmr}: given a tableau $\ttt$ and $1\ls i,j\ls n$, we write $i\to_\ttt j$ to mean that $i$ and $j$ lie in the same row of the same component, with $j$ to the right of $i$. We write $i\nearrow_\ttt j$ to mean that $i$ and $j$ lie in the same component of $\ttt$, with $j$ strictly higher and strictly to the right, and we write $i\Nearrow_\ttt j$ to mean that either $i\nearrow_\ttt j$ or $j$ lies in an earlier component than $i$.  The notations $i\downarrow_\ttt j$, $i\swarrow_\ttt j$ and $i\Swarrow_\ttt j$ are defined similarly.

There are two standard $\la$-tableaux of particular importance. The tableau $\ttt_\la$ is the standard tableau obtained by writing $1,\dots,n$ in order down successive columns from left to right, while $\ttt^\la$ is the tableau obtained by writing $1,\dots,n$ in order along successive rows from top to bottom. Note that we then have $\ttt^\la=(\ttt_{\la'})'$.

\begin{eg}
With $\la=\left((2^2),(2,1^2),(3,2)\right)$ we have
\[
\ttt_\la=\gyoung(^6;\ten;\twelve,^6;\eleven;\thirteen,/\hf,^3^\hf;6;9,^3^\hf;7,^3^\hf;8,/\hf,;1;3;5,;2;4),\qquad
\ttt^\la=\gyoung(^6;1;2,^6;3;4,/\hf,^3^\hf;5;6,^3^\hf;7,^3^\hf;8,/\hf,;9;\ten;\eleven,;\twelve;\thirteen).
\]
\end{eg}

The symmetric group $\sss n$ acts naturally on the left on the set of $\la$-tableaux. Given a $\la$-tableau $\ttt$, we define the permutations $w_\ttt$ and $w^\ttt$ in $\sss n$ by
\[
w_\ttt \ttt_\la = \ttt = w^\ttt \ttt^\la.
\]

Later we shall need the following lemma; recall that $\ls\Lf$ denotes the left order on $\sss n$.

\begin{lemma}\label{leftstd}
Suppose $\la\in\mptn ln$ and $\tts,\ttt$ are $\la$-tableaux with $w_\tts\ls\Lf w_\ttt$. If $\ttt$ is standard, then $\tts$ is standard.
\end{lemma}

\begin{pf}
Using induction on $l(w_\ttt)-l(w_\tts)$, we may assume $l(w_\ttt)=l(w_\tts)+1$, which means in particular that $\ttt=s_i\tts$ for some $i$. Since $\ttt$ is standard, the only way $\tts$ could fail to be standard is if $i+1$ occupies the node immediately below or immediately to the right of $i$ in $\ttt$. But either possibility means that $i$ occurs before $i+1$ in the `column reading word' of $\ttt$, i.e.\ the word obtained by reading the entries of $\ttt$ down successive columns from left to right. In other words, $w_\ttt^{-1}(i)<w_\ttt^{-1}(i+1)$, but this means that $l(w_\tts)>l(w_\ttt)$, a contradiction.
\end{pf}

Now we introduce a dominance order on tableaux. If $\tts,\ttt$ are $\la$-tableaux, then we write $\tts\dom\ttt$ if and only if $w_\tts\cgs w_\ttt$ (recall that $\cgs$ denotes the Bruhat order on $\sss n$). There should be no ambiguity in using the symbol $\dom$ for both the dominance order on multipartitions and the dominance order on tableaux.

There is an alternative description of the dominance order on tableaux which will be very useful. If $\ttt$ is a $\la$-tableau and $0\ls m\ls n$, we define $\tabupto\ttt m$ to be the set of nodes of $[\la]$ whose entries are less than or equal to $m$. If $\ttt$ is row-strict, then $\tabupto\ttt m$ is the Young diagram of an $l$-multicomposition of $m$, which we call $\shp\ttt m$. If $\ttt$ is standard, then $\shp\ttt m$ is an $l$-multipartition of $m$.

Now we have the following \lcnamecref{brudom}. This is proved in the case $l=1$ in \cite[Theorem 3.8]{math} (where it is attributed to Ehresmann and James); in fact, the proof in \cite{math} carries over to the case of arbitrary $l$ without any modification.

\begin{propn}\label{brudom}
Suppose $\la\in\mptn ln$ and $\tts,\ttt$ are row-strict $\la$-tableaux. Then $\tts\domby\ttt$ if and only if $\shp\tts m\domby\shp\ttt m$ for $m=1,\dots,n$.
\end{propn}

In this paper, we shall briefly consider a natural analogue of this notion for column-strict tableaux. Suppose $\la\in\mptn ln$ and $\ttt$ is a column-strict $\la$-tableau; define the diagram $\tabupto\ttt m$ as above, and define $\tabupto\ttt m'$ to be the `conjugate diagram' to $\tabupto\ttt m$, that is
\[
\tabupto\ttt m'=\rset{(c,r,l+1-k)}{(r,c,k)\in\tabupto\ttt m}.
\]
Then $\tabupto\ttt m'$ is the Young diagram of an $l$-multicomposition of $m$, which we denote $\shp\ttt m'$. Now we have following statement, which can be deduced from \cref{brudom} by conjugating tableaux.

\begin{propn}\label{brudomconj}
Suppose $\la\in\mptn ln$ and $\tts,\ttt$ are column-strict $\la$-tableaux. Then $\tts\domby\ttt$ if and only if $\shp\tts m'\dom\shp\ttt m'$ for $m=1,\dots,n$.
\end{propn}


\subsection{Residues and degrees}\label{resdegsec}

In this section we connect the Lie-theoretic set-up above with multipartitions and tableaux. We fix an $e$-multicharge $\kappa=(\kappa_1,\dots,\kappa_l)$. We define the \emph{residue} $\res A=\res^\kappa A$ of a node $A=(r,c,m)\in\bbn\times\bbn\times\{1,\dots,l\}$ by
\[
\res A = \kappa_m + (c-r)\pmod e.
\]
We say that $A$ is an \emph{$i$-node} if it has residue $i$. Given $\la\in\mptn ln$, we define the \emph{content} of $\la$ to be the root
\[
\cont(\la)=\sum_{A\in[\la]}\alpha_{\res A}.
\]
We then define the \emph{defect} $\df(\la)$ of $\la$ to be $\df(\cont(\la))$.

If $\ttt$ is a $\la$-tableau, we define its \emph{residue sequence} to be the sequence $i(\ttt)=(i_1,\dots,i_n)$, where $i_r$ is the residue of the node $\ttt^{-1}(r)$, for each $r$. The residue sequences of the tableaux $\ttt_\la$ and $\ttt^\la$ will be of particular importance, and we set $i_\la:=i(\ttt_\la)$ and $i^\la:=i(\ttt^\la)$.

\begin{eg}
Take $\la=\left((2^2),(2,1^2),(3,2)\right)$ as in the last example, and suppose $e=4$ and $\kappa=(1,2,0)$. Then the residues of the nodes of $\la$ are given by the following diagram.
\[
\gyoung(^6;1;2,^6;0;1,/\hf,^3^\hf;2;3,^3^\hf;1,^3^\hf;0,/\hf,;0;1;2,;3;0)
\]
So we have
\[
i_\la=(0,3,1,0,2,2,1,0,3,1,0,2,1),\qquad i^\la=(1,2,0,1,2,3,1,0,0,1,2,3,0).
\]
\end{eg}

Now we recall from \cite[\S3.5]{bkw} the degree and codegree of a standard tableau. Suppose $\la\in\mptn ln$ and $A$ is an $i$-node of $\la$. Set
\[
d_A(\la):=\big|\{\text{addable $i$-nodes of $\la$ strictly below }A\}\big|-\big|\{\text{removable $i$-nodes of $\la$ strictly below }A\}\big|,
\]
and
\[
d^A(\la):=\big|\{\text{addable $i$-nodes of $\la$ strictly above }A\}\big|-\big|\{\text{removable $i$-nodes of $\la$ strictly above }A\}\big|.
\]

For $\ttt\in\std\la$ we define the \emph{degree} of $\ttt$ recursively, setting $\deg(\ttt):=0$ when $\ttt$ is the unique $\varnothing$-tableau. If $\ttt\in\std\la$ with $|\la|>0$, let $A=\ttt^{-1}(n)$, let $\ttt_{<n}$ be the tableau obtained by removing this node and set
\[
\deg(\ttt):= d_A(\la)+\deg(\ttt_{<n}).
\]
Similarly, define the \emph{codegree} of $\ttt$ by setting $\codeg(\ttt):=0$ if $\ttt$ is the unique $\varnothing$-tableau, and
\[
\codeg(\ttt):= d^A(\la)+\codeg(T_{<n})
\]
for $\ttt\in\std\la$ with $|\la|>0$. We note that the definitions of degree and codegree depend on the \emc $\kappa$, and therefore we write $\deg^\kappa$ and $\codeg^\kappa$ when we wish to emphasise $\kappa$.

\begin{eg}
Suppose $e=3$, $\kappa=(1,1)$ and $\ttt$ is the $\left((2),(2,1)\right)$-tableau
\[
\gyoung(^2^\hf;3;4,/\hf,;1;5,;2)
\]
which has residue sequence $i(\ttt)=(1,0,1,2,2)$. Letting $A=\ttt^{-1}(5)=(1,2,2)$, we find that $d_A(\la)=1$ and $d^A(\la)=-1$. Recursively one finds that for the tableau
\[
\ttt_{<5}=\gyoung(^2^\hf;3;4,/\hf,;1,;2)
\]
we have $\deg(\ttt_{<5})=2$ and $\codeg(\ttt_{<5})=1$, so that $\deg(\ttt)=3$ and $\codeg(\ttt)=0$.
\end{eg}

The degree, codegree of a standard $\la$-tableau are related to the defect of $\la$ by the following result.

\begin{lemmac}{bkw}{Lemma 3.12}\label{degdef}
Suppose $\la\in\mptn ln$ and $\ttt\in\std\la$. Then
\[
\deg(\ttt)+\codeg(\ttt)=\df(\la).
\]
\end{lemmac}

\subsection{KLR algebras}\label{klrsec}

We now give the definition of the algebras which will be our main object of study.

Suppose $\alpha$ is a positive root of height $n$, and set
\[
I^\alpha=\lset{i\in I^n}{\alpha_{i_1}+\dots+\alpha_{i_n}=\alpha}.
\]
Now define $\hhh\alpha$ to be the unital associative $\bbf$-algebra with generating set
\[
\lset{e(i)}{i\in I^\alpha}\cup\{y_1,\dots,y_n\}\cup\{\psi_1,\dots,\psi_{n-1}\}
\]
and relations
{\allowdisplaybreaks
\begin{align*}
e(i)e(j)&=\delta_{i,j} e(i);\\
\sum_{i \in I^\alpha} e(i)&=1;\\
y_re(i)&=e(i)y_r;\\
\psi_r e(i) &= e(s_ri) \psi_r;\\
y_ry_s&=y_sy_r;\\
\psi_ry_s&=\mathrlap{y_s\psi_r}\hphantom{\smash{\begin{cases}(\psi_{r+1}\psi_r\psi_{r+1}+y_r-2y_{r+1}+y_{r+2})e(i)\\\\\\\end{cases}}}\kern-\nulldelimiterspace\text{if } s\neq r,r+1;\\
\psi_r\psi_s&=\mathrlap{\psi_s\psi_r}\hphantom{\smash{\begin{cases}(\psi_{r+1}\psi_r\psi_{r+1}+y_r-2y_{r+1}+y_{r+2})e(i)\\\\\\\end{cases}}}\kern-\nulldelimiterspace\text{if } |r-s|>1;\\
y_r \psi_r e(i) &=(\psi_r y_{r+1} - \delta_{i_r,i_{r+1}})e(i);\\
y_{r+1} \psi_r e(i) &=(\psi_r y_r + \delta_{i_r,i_{r+1}})e(i);\\
\psi_r^2 e(i)&=\begin{cases}
\mathrlap0\phantom{(\psi_{r+1}\psi_r\psi_{r+1}+y_r-2y_{r+1}+y_{r+2})e(i)}& \text{if }i_r=i_{r+1},\\
e(i) & \text{if }i_{r+1}\neq i_r, i_r\pm1,\\
(y_{r+1} - y_r) e(i) & \text{if }i_r \rightarrow i_{r+1},\\
(y_r - y_{r+1}) e(i) & \text{if }i_r\leftarrow i_{r+1},\\
(y_{r+1} - y_r)(y_r - y_{r+1}) e(i) & \text{if }i_r\rightleftarrows i_{r+1};
\end{cases}\\
\psi_r\psi_{r+1}\psi_re(i)&=\begin{cases}
(\psi_{r+1}\psi_r\psi_{r+1}+1)e(i)& \text{if }i_{r+2}=i_r\rightarrow i_{r+1},\\
(\psi_{r+1}\psi_r\psi_{r+1}-1)e(i)& \text{if }i_{r+2}=i_r\leftarrow i_{r+1},\\
(\psi_{r+1}\psi_r\psi_{r+1}+y_r-2y_{r+1}+y_{r+2})e(i)& \text{if }i_{r+2}=i_r\rightleftarrows i_{r+1},\\
(\psi_{r+1}\psi_r\psi_{r+1})e(i)& \text{otherwise;}
\end{cases}
\end{align*}
}for all admissible $r,s,i,j$.

The \emph{affine Khovanov--Lauda--Rouquier algebra} or \emph{quiver Hecke algebra} $\hhh n$ is defined to be the direct sum $\bigoplus_\alpha\hhh\alpha$, where the sum is taken over all positive roots of height $n$.

\begin{rmks}\indent
\begin{enumerate}
\vspace{-\topsep}
\item
We use the same notation for the generators $\psi_r$ and $y_s$ for different $\alpha$; when using these generators, we shall always make it clear which algebra $\hhh\alpha$ these generators are taken from.
\item
When $e<\infty$, we can modify the above presentation of $\hhh\alpha$ to give a presentation for $\hhh n$: we take generating set $\lset{e(i)}{i\in I^n}\cup\{y_1,\dots,y_n\}\cup\{\psi_1,\dots,\psi_{n-1}\}$, and replace the relation $\sum_{i\in I^\alpha}e(i)=1$ with $\sum_{i\in I^n}e(i)=1$. The generator $\psi_r$ in this presentation is just the sum of the corresponding generators $\psi_r$ of the individual algebras $\hhh\alpha$ in the direct sum $\bigoplus_\alpha\hhh\alpha$, and similarly for $y_s$. When $e=\infty$ we cannot do this, since the set $I^n$ is infinite (in fact, $\hhh n$ is non-unital in this case).
\end{enumerate}
\end{rmks}

The following result can easily be checked from the definition of $\hhh\alpha$.

\begin{lemmac}{bk}{Corollary 1}
There is a $\bbz$-grading on the algebra $\hhh\alpha$ such that for all admissible $r$ and $i$,
\[
\deg(e(i))=0, \quad \deg(y_r)=2, \quad \deg(\psi_re(i))=-a_{i_ri_{r+1}}.
\]
\end{lemmac}

\subsubsection*{Shift maps}

Recall from \S\ref{symsec} that $\shift k:\sss m\to\sss n$ denotes the homomorphism defined by $s_i\mapsto s_{i+k}$. We now define the corresponding maps for the algebras $\hhh\alpha$.

\begin{defn}
Suppose $1\leq m \leq n$ and $0\leq k\leq n-m$, and that $\alpha$ and $\beta$ are positive roots with $\alpha$ of height $n$ and $\beta$ of height $m$. Given $i\in I^\beta$, define $J_i:=\lset{j\in I^\alpha}{j_{s+k}=i_s\text{ for }1\leq s \leq m}$, and let $e(i)^{+k}=\sum_{j\in J_i}e(j)$. Now define the homomorphism $\shift k:\hhh\beta\rightarrow\hhh\alpha$ by
\[
e(i)\mapsto e(i)^{+k},\qquad\psi_re(i)\mapsto\psi_{r+k}e(i)^{+k},\qquad y_re(i)\mapsto y_{r+k}e(i)^{+k}.
\]
\end{defn}
It is easy to check from the definition of $\hhh\alpha$ that $\shift k$ is a degree-preserving (non-unital) homomorphism of algebras. Moreover, the PBW-type basis theorem for $\hhh\alpha$ in \cite[Theorem 2.5]{kl} and \cite[Theorem 3.7]{rouq} shows that if $\beta\ls\alpha$ then $\shift k$ is injective (obviously $\shift k$ is the zero map if $\beta\nleqslant\alpha$).

\subsubsection*{Cyclotomic algebras and the Brundan--Kleshchev isomorphism theorem}

Given a positive root $\alpha$ and an \emc $\kappa=(\kappa_1,\dots,\kappa_l)\in I^l$, we define $\hhh\alpha^\kappa$ to be the quotient of $\hhh\alpha$ by the \emph{cyclotomic relations}
\[
y_1^{\bil{\La_\kappa}{\alpha_{i_1}}}e(i)=0\quad\text{for } i\in I^\alpha.
\]
The \emph{cyclotomic KLR algebra} $\hhh n^\kappa$ is then defined to be the sum $\bigoplus_\alpha\hhh\alpha^\kappa$. Here we sum over all positive roots $\alpha$ of height $n$, though in fact only finitely many of the summands will be non-zero, so (even when $e=\infty$) $\hhh n^\kappa$ is a unital algebra.

Note that the embedding $\shift 0$ passes naturally into the cyclotomic quotients.

A stunning result of Brundan and Kleshchev \cite[Main Theorem]{bk} is that if $e=\infty$ or $e$ is not divisible by $\nchar(\bbf)$, then $\hhh n^\kappa$ is isomorphic to an Ariki--Koike algebra of level $l$, defined at an $e$th root of unity. Similarly, if $e=\nchar(\bbf)$, then $\hhh n^\kappa$ is isomorphic to a degenerate Ariki--Koike algebra; in particular, when $l=1$, $\hhh n^\kappa$ is isomorphic to the group algebra $\bbf\sss n$. As a consequence, these Hecke algebras are non-trivially $\bbz$-graded. This theorem motivates our choice of notation $\hhh n$ for the KLR algebra.

\subsection{Specht modules}\label{spechtsec}

We now recall the universal graded row and column Specht modules introduced by Kleshchev, Mathas and Ram; we refer the reader to \cite[\S\S5,7]{kmr} for further details.

Fix an \emc $\kappa$. Suppose $\la\in\mptn ln$, and let $\alpha=\cont(\la)$. Say that a node $A=(r,c,m)\in[\la]$ is a \emph{column Garnir node} if $(r,c+1,m)\in[\la]$. The \emph{column Garnir belt} $\bfB_A$ is defined to be the set of nodes
\[
\bfB_A=\lset{(s,c,m)\in[\la]}{s\geq r}\cup\lset{(s,c+1,m)\in[\la]}{s\leq r}.
\]
This belt is used to define a \emph{column Garnir element} $\garn A\in\hhh\alpha$. The full definition of $\garn A$ is quite complicated, and can be found in \cite[Definition 7.10]{kmr}. Here we just give $\garn A$ explicitly in a special case which we will use in the proof of \cref{wjdom}, and record some useful properties of $\garn A$ which apply in general.

For our special case, we suppose that $A$ is a Garnir node of $\la$ of the form $(1,c,m)$. If $a$ is the entry in node $A$ of $\ttt_\la$ and $b$ is the entry in node $(1,c+1,m)$, then $\garn A=\psi_a\psi_{a+1}\dots\psi_{b-1}$.

Now suppose $A=(r,c,m)$ is an arbitrary Garnir node of $\la$. Then in $\ttt_\la$ the nodes of $\bfB_A$ are occupied by the integers $a,a+1,\dots,b$ for some $a<b$. The following facts can be distilled from \cite[\S7]{kmr}:
\begin{itemize}
\item
$\garn A$ is a linear combination of products of the form $\psi_{i_1}\dots\psi_{i_d}$ where $a\leq i_1,\dots,i_d<b$;
\item
$\garn A$ depends only on $e$, $r$, $a$ and the length of the column containing $A$.
\end{itemize}
(In fact, as defined in \cite{kmr}, the column Garnir element $\garn A$ also involves an idempotent $e(i)$ which depends on $\la$ and makes $\garn A$ homogeneous, but this term can be omitted without affecting the Garnir relation given below.)

\begin{eg}
For example, let $\la=((3,3,2,2,1),(2,1))$ and let $A=(3,1,1)$. Then $A$ is a column Garnir node, and $\ttt_\la$ (with the Garnir belt $\bfB_A$ shaded) is as follows.
\Yfillcolour{white!80!black}
\[
\ttt_\la=\yfz
\raisebox{1.2cm}{
\gyoung(^\tpf;4!\yfo;9!\yfz;<13>,^\tpf;5!\yfo;<10>!\yfz;<14>,^\tpf!\yfo;6;\eleven,^\tpf;7!\yfz;<12>,^\tpf!\yfo;8,!\yfz/\hf,;1;3,;2)
}.
\]
The column Garnir element $\garn A$ is then a linear combination of products of the generators $\psi_6,\psi_7,\psi_8,\psi_9,\psi_{10}$; the exact expression for $\garn A$ depends on the choice of $e$.
\end{eg}

Now define the \emph{column Specht module} $\spek\la\kappa$ to be the graded $\hhh\alpha$-module generated by the vector $z_\la$ of degree $\codeg(\ttt_\la)$ subject to the following relations:

\begin{enumerate}
\item $e(i_\la)z_\la = z_\la$;
\item $y_r z_\la = 0$ for all $r=1,\dots,n$;
\item $\psi_r z_\la = 0$ for all $r=1,\dots,n-1$ such that $r\downarrow_{\ttt_\la}r+1$;
\item $\garn A z_\la=0$ for all column Garnir nodes $A\in\la$.
\end{enumerate}

We may relax notation and just write $\spe\la$, if the \emc $\kappa$ is understood. We shall mostly consider $\spe\la$ as an $\hhh n$-module, by setting $\hhh\beta\spe\la=0$ for $\beta\neq\alpha$. Thus we have $\hhh n$-modules $\spek\la\kappa$ for all \emcs $\kappa$ and all $\la\in\mptn ln$. The main purpose of this paper is to study the space of $\hhh n$-homomorphisms $\spe\la\to\spe\mu$, for $\la,\mu\in\mptn ln$. The following result is obvious from the definitions.

\begin{lemma}\label{samedef}
Suppose $\la,\mu\in\mptn ln$, and let $\alpha=\cont(\la)$. If $\shom\la\mu\neq0$, then $\cont(\mu)=\alpha$ (and in particular $\df(\la)=\df(\mu)$), and $\shom\la\mu=\hom_{\hhh\alpha}(\spe\la,\spe\mu)$.
\end{lemma}


We shall also need to consider row Specht modules; for these, the definitions are largely obtained by `conjugating' the definitions for column Specht modules. Fix $\kappa$, $\la$ and $\alpha$ as above. Say that a node $A=(r,c,m)\in[\la]$ is a \emph{row Garnir node} if $(r+1,c,m)\in[\la]$, and define the \emph{row Garnir belt}
\[
\bfB^A=\lset{(r,d,m)\in[\la]}{d\gs b}\cup\lset{(r+1,d,m)\in[\la]}{d\ls c}.
\]
This belt is used to define a \emph{row Garnir element} $\rgarn A$. We refer the reader to \cite[Definition 5.8]{kmr} for the definition of this; here we just note the following facts:
\begin{itemize}
\item
in $\ttt^\la$ the nodes of $\bfB^A$ are occupied by the integers $a,a+1,\dots,b$ for some $a<b$;
\item
$\rgarn A$ is a linear combination of products of the form $\psi_{i_1}\dots\psi_{i_d}$ where $a\leq i_1,\dots,i_d<b$;
\item
$\rgarn A$ depends only on $e$, $c$, $a$ and the length of the row containing $A$.
\end{itemize}

Now we can define the \emph{row Specht module} $\rspe\la$, which is the graded $\hhh\alpha$-module generated by the vector $z^\la$ of degree $\deg(\ttt^\la)$ subject to the relations

\begin{enumerate}
\item $e(i^\la)z^\la = z^\la$;
\item $y_rz^\la = 0$ for all $r=1,\dots,n$;
\item $\psi_rz^\la = 0$ for all $r=1,\dots,n-1$ such that $r\rightarrow_{\ttt^\la}r+1$;
\item $\rgarn Az^\la=0$ for all row Garnir nodes $A\in\la$.
\end{enumerate}

We define basis elements for the row and column Specht modules as follows. For each $\ttt\in\std\la$ we fix a \emph{preferred reduced expression} $s_{r_1}\dots s_{r_a}$ for the permutation $w_\ttt$, and define $\psi_\ttt:=\psi_{r_1}\dots\psi_{r_a}$ and $v_\ttt:=\psi_\ttt z_\la$. Similarly, we fix a preferred reduced expression $s_{t_1}\dots s_{t_b}$ for $w^\ttt$, and set $\psi^\ttt:=\psi_{t_1}\dots\psi_{t_b}$ and $v^\ttt:=\psi^\ttt z^\la$.

Note that the elements $v_\ttt$ and $v^\ttt$ may depend on the choice of preferred reduced expressions, since the $\psi_r$ do not satisfy the braid relations. However, the following results are independent of the choices made.

\begin{lemmac}{kmr}{Propositions 5.14 \& 7.14}
Suppose $\la\in\mptn ln$ and $\ttt\in\std\la$. Then $\deg(v^\ttt)=\deg(\ttt)$ and $\deg(v_\ttt)=\codeg(\ttt)$.
\end{lemmac}

\begin{lemmac}{kmr}{Corollaries 6.24 \& 7.20}
Suppose $\la\in\mptn ln$. Then $\rset{v^\ttt}{\ttt\in\std\la}$ is an $\bbf$-basis for $\rspe\la$, and $\rset{v_\ttt}{\ttt\in\std\la}$ is an $\bbf$-basis for $\spe\la$.
\end{lemmac}

In spite of the dependence of these bases on the choices of preferred reduced expressions, we refer to the bases $\rset{v^\ttt}{\ttt\in\std\la}$ and $\rset{v_\ttt}{\ttt\in\std\la}$ as the \emph{standard bases} for $\rspe\la$ and $\spe\la$ respectively.

For the remainder of this section we summarise some basic results about the action of $\hhh\alpha$ on $\spe\la$. Many of these results are cited from \cite{bkw}, where they are stated for row Specht modules. In this paper we concentrate as far as possible on column Specht modules, so we translate all the results to this setting. Throughout we fix $\la\in\mptn ln$, and let $\psi_1,\dots,\psi_{n-1}$ refer to the generators of $\hhh\alpha$, where $\alpha=\cont(\la)$. Recall that if $\tts,\ttt$ are standard $\la$-tableaux, then we write $\tts\dom\ttt$ to mean that $w_\tts\cgs w_\ttt$.

\begin{lemmac}{bkw}{Theorem 4.10(i)}\label{diffredex}
Suppose $\ttt\in\std\la$, and $s_{j_1}\dots s_{j_r}$ is any reduced expression for $w_\ttt$. Then $\psi_{j_1}\dots\psi_{j_r}z_\mu-v_\ttt$ is a linear combination of basis elements $v_\ttu$ for $\ttu\domsby\ttt$.
\end{lemmac}

\begin{lemmac}{bkw}{Lemma 4.9}\label{samecol}
Suppose $\ttt\in\std\la$ and that $j-1\to_\ttt j$ or $j-1\downarrow_\ttt j$. Then $\psi_{j-1}v_\ttt$ is a linear combination of basis elements $v_\ttu$ for $\ttu\domsby\ttt$.
\end{lemmac}

\begin{lemmac}{bkw}{Lemma 4.8}\label{ystd}
Suppose $\ttt\in\std\la$ and $1\ls i\ls n$. Then $y_iv_\ttt$ is a linear combination of basis elements $v_\ttu$ for $\ttu\domsby\ttt$.
\end{lemmac}

We'll use \cref{diffredex,ystd} to prove the following similar result, which is suggested but not proved in the proof of \cite[Theorem 4.10]{bkw}.

\begin{lemma}\label{unswitch}
Suppose $\ttt\in\std\la$ and $j-1\Swarrow_\ttt j$. Then $\psi_{j-1}v_\ttt$ is a linear combination of basis elements $v_\ttu$ for $\ttu\domby\ttt$.
\end{lemma}

We begin with the following simple observation.

\begin{lemma}\label{jleft}
Suppose $\ttt\in\std\la$. Then $j-1\Swarrow_\ttt j$ if and only if $w_\ttt$ has a reduced expression beginning with $s_{j-1}$.
\end{lemma}

\begin{pf}
Both conditions are equivalent to the condition that $w_\ttt^{-1}(j-1)>w_\ttt^{-1}(j)$.
\end{pf}

\begin{pf}[Proof of \cref{unswitch}]
By \cref{jleft}, $w_\ttt$ has a reduced expression of the form $s_{j-1}s_{k_1}\dots s_{k_r}$. Using \cref{diffredex} we have
\[
v_\ttt=\psi_{j-1}\psi_{k_1}\dots\psi_{k_r}z_\la+\sum_{\substack{\ttu\in\std\la\\\ttu\domsby\ttt}}a_\ttu v_\ttu
\]
for some $a_\ttu\in\bbf$. So
\[
\psi_{j-1}v_\ttt=\psi_{j-1}^2\psi_{k_1}\dots\psi_{k_r}z_\la+\sum_{\substack{\ttu\in\std\la,\\\ttu\domsby\ttt}}a_\ttu\psi_{j-1}v_\ttu.\tag*{\cnd}
\]
Using the KLR relations (and moving the appropriate idempotent $e(i)$ through), the first term on the right-hand side becomes $g\psi_{k_1}\dots\psi_{k_r}z_\la$, where $g$ is a polynomial in $y_1,\dots,y_n$. Now $s_{k_1}\dots s_{k_r}$ is a reduced expression for the standard tableau $\tts=s_{j-1}\ttt$, so by \cref{diffredex} we have
\[
\psi_{k_1}\dots\psi_{k_r}z_\la=v_\tts+\sum_{\substack{\ttv\in\std\la,\\\ttv\domsby\tts}}b_\ttv v_\ttv
\]
for some $b_\ttv\in\bbf$. So (since $\tts\domsby\ttt$) the first term on the right-hand side of \cnd{} is a linear combination of terms of the form $g v_\ttv$ for $\ttv\in\std\la$ with $\ttv\domsby\ttt$. By \cref{ystd} this reduces to a linear combination of basis elements $v_\ttv$ for $\ttv\domsby\ttt$.

Now consider each of the remaining terms $\psi_{j-1}v_{\ttu}$ in \cnd. If $j-1\Swarrow_\ttu j$, then by induction on the Bruhat order $\psi_{j-1}v_{\ttu}$ is a linear combination of basis elements $v_\ttv$ for $\ttv\domby\ttu\domsby\ttt$, so we can ignore any such $\ttu$. If $j-1\to_\ttu j$ or $j-1\downarrow_\ttu j$, then we apply \cref{samecol} to get the same conclusion. If $j-1\Nearrow_\ttu j$, let $\ttr$ be the tableau obtained by swapping $j-1$ and $j$ in $\ttu$; then a reduced expression for $w_\ttr$ may be obtained by adding $s_{j-1}$ at the start of a reduced expression for $w_\ttu$, and we have $\ttr\domby\ttt$ by \cref{brulem}. So by \cref{diffredex} again,
\[
\psi_{j-1}v_\ttu = v_\ttr + \sum_{\ttw\domsby\ttr}c_\ttw v_\ttw
\]
for some $c_\ttw\in\bbf$, and we are done.
\end{pf}

\begin{lemma}\label{subexp}
Suppose $\la\in\mptn ln$, and $\ttt\in\std\la$. Suppose $j_1,\dots,j_r\in\{1,\dots,n-1\}$, and that when $\psi_{j_1}\dots\psi_{j_r}z_\la$ is expressed as a linear combination of standard basis elements, $v_\ttt$ appears with non-zero coefficient. Then the expression $s_{j_1}\dots s_{j_r}$ has a reduced expression for $w_\ttt$ as a subexpression.
\end{lemma}

\begin{pf}
We proceed by induction on $r$, with the case $r=0$ trivial. Let $j=j_1$. Then by assumption $v_\ttt$ appears with non-zero coefficient in $\psi_jv_\tts$, where $\tts\in\std\la$ and $v_\tts$ appears with non-zero coefficient in $\psi_{j_2}\dots\psi_{j_r}z_\la$. By induction the expression $s_{j_2}\dots s_{j_r}$ has a subexpression which is a reduced expression for $w_\tts$, so if $w_\ttt\ls w_\tts$ (i.e.\ if $\ttt\domby\tts$) then we are done. By \cref{samecol,unswitch}, this happens if $j\rightarrow_\tts j+1$, $j\downarrow_\tts j+1$ or $j\Swarrow_\tts j+1$. So we can assume that $j\Nearrow_\tts j+1$. But in this case $w_\ttt=s_jw_\tts$, with $l(w_\ttt)=l(w_\tts)+1$, so $w_\ttt$ has a reduced expression obtained by adding $s_j$ at the start of a reduced expression for $w_\tts$. So again the result follows by induction.
\end{pf}

\subsection{Specht modules for $\hhh n^\kappa$ and homomorphisms}

Throughout this paper we consider the Specht module $\spe\la$ as a module for the \emph{affine} algebra $\hhh\alpha$ (where $\alpha=\cont(\la)$) and by extension for the algebra $\hhh n$. In fact, it is not hard to show that $\spe\la$ is annihilated by the element $y_1^{\bil{\La_\kappa}{\alpha_{i_1}}}e(i)$ for every $i$, so that $\spe\la$ is a module for the cyclotomic algebra $\hhh n^\kappa$ introduced in \S\ref{klrsec}. We shall almost entirely be studying the space of $\hhh n$-homomorphisms between two Specht modules $\spe\la$ and $\spe\mu$ defined for the same \emc $\kappa$, and clearly in this situation $\hhh n$-homomorphisms between these two modules are the same as $\hhh n^\kappa$-homomorphisms. In view of the Brundan--Kleshchev isomorphism theorem mentioned above, our results can therefore be viewed as statements about homomorphisms between Specht modules for (degenerate) Ariki--Koike algebras, and so they generalise the results of the first author and Lyle for homomorphisms between Specht modules for the symmetric group \cite[Theorem 2.1]{fl}, and of Lyle and Mathas for Hecke algebras of type $A$ \cite[Theorem 1.1]{lm}.

In this paper, however, we restrict attention entirely to the affine algebra $\hhh n$. This is because we occasionally (in particular, in \cref{sign}) need to compare Specht modules defined for different \emcs.

\section{$\la$-dominated tableaux and dominated homomorphisms}\label{domsec}

In this paper we consider the space of homomorphisms between two given Specht modules. However, our results concerning row and column removal will only apply to homomorphisms of a certain kind, which we call \emph{dominated} homomorphisms. But as we shall see in \cref{domhom}, in many cases all homomorphisms between Specht modules are dominated.

\subsection{$\la$-dominated tableaux}\label{domtsec}

Suppose $\la,\mu\in\mptn ln$ and $\ttt\in\std\mu$. Given $0\ls j\ls n$, we say that $\ttt$ is \emph{$\la$-column-dominated on $1,\dots,j$} if each $i\in\{1,\dots,j\}$ appears at least as far to the left in $\ttt$ as it does in $\ttt_\la$. We say simply that $\ttt$ is \emph{$\la$-column-dominated} if it is $\la$-column-dominated on $1,\dots,n$. We remind the reader of our unusual convention for drawing Young diagrams, in which a node $(r,c,m)$ lies to the left of $(r',c',m')$ if either $m>m'$ or ($m=m'$ and $c\ls c'$).

We write $\cstd\la\mu$ for the set of $\la$-column-dominated standard $\mu$-tableaux. It is easy to see that $\cstd\la\mu$ is non-empty if and only if $\la\dom\mu$, and that $\cstd\mu\mu=\{\ttt_\mu\}$.

We say $\ttt$ is \emph{weakly $\la$-column-dominated on $1,\dots,j$} if each $i\in\{1,\dots,j\}$ appears in a component at least as far to the left in $\ttt$ as it does in $\ttt_\la$. We say that $\ttt$ is weakly $\la$-column-dominated if it is weakly $\la$-column-dominated on $1,\dots,n$.

We also introduce row-dominance. Say that $\ttt\in\std\mu$ is \emph{$\la$-row-dominated} if each $i\in\{1,\dots,n\}$ appears at least as high in $\ttt$ as it does in $\ttt^\la$. We write $\rstd\la\mu$ for the set of $\la$-row-dominated standard $\mu$-tableaux, which is non-empty if and only if $\la\domby\mu$.

Since we shall primarily be considering column Specht modules, we shall often simply say `$\la$-dominated' meaning `$\la$-column-dominated'.

We give a helpful alternative characterisation of the $\la$-dominated and $\la$-row-dominated properties.

\begin{lemma}\label{domalt}
Suppose $\la,\mu\in\mptn ln$, and $\tts\in\std\mu$.
\begin{enumerate}
\item\label{domalt1}
$\tts$ is $\la$-column-dominated on $1,\dots,j$ if and only if $\shp{(\ttt_\la)}m\dom\shp\tts m$ for all $m=1,\dots,j$.
\item\label{domalt2}
$\tts$ is $\la$-row-dominated on $1,\dots,j$ if and only if $\shp{(\ttt^\la)}m\domby\shp\tts m$ for all $m=1,\dots,j$.
\end{enumerate}
\end{lemma}

\begin{pf}
We prove only (\ref{domalt2}); the proof of (\ref{domalt1}) is analogous. Suppose first that $\tts$ is not $\la$-row-dominated on $1,\dots,j$. Choose an entry $m\leq j$ which appears strictly lower in $\tts$ than in $\ttt^\la$, and let $\tau=\shp{(\ttt^\la)}m$ and $\sigma=\shp\tts m$. Suppose that $m$ appears in position $(r,c,k)$ in $\ttt^\la$. The construction of $\ttt^\la$ means that the entries $1,\dots,m-1$ all appear at least as high as $m$ in $\ttt^\la$, and so
\[
|\tau^{(1)}|+\dots+|\tau^{(k-1)}|+\tau^{(k)}_1+\dots+\tau^{(k)}_r=m.
\]
On the other hand, $m$ appears below row $r$ of component $k$ in $\tts$, so
\[
|\sigma^{(1)}|+\dots+|\sigma^{(k-1)}|+\sigma^{(k)}_1+\dots+\sigma^{(k)}_r<m.
\]
Hence $\tau\ndomby\sigma$.

Conversely, suppose $\shp{\ttt^\la}m\ndomby\shp\tts m$ for some $m\leq j$; choose such an $m$, and let $\tau=\shp{\ttt^\la}m$ and $\sigma=\shp\tts m$. Since $\tau\ndomby\sigma$, there are $r,k$ such that
\[
|\tau^{(1)}|+\dots+|\tau^{(k-1)}|+\tau^{(k)}_1+\dots+\tau^{(k)}_r>|\sigma^{(1)}|+\dots+|\sigma^{(k-1)}|+\sigma^{(k)}_1+\dots+\sigma^{(k)}_r.
\]
If we let $d=|\tau^{(1)}|+\dots+|\tau^{(k-1)}|+\tau^{(k)}_1+\dots+\tau^{(k)}_r$, then $d\ls m$ and the integers $1,\dots,d$ all appear in row $r$ of component $k$ or higher in $\ttt^\la$. Since $|\sigma^{(1)}|+\dots+|\sigma^{(k-1)}|+\sigma^{(k)}_1+\dots+\sigma^{(k)}_r<d$, at least one of the integers $1,\dots,d$ appears in $\tts$ below row $r$ of component $k$. So there is some $i\ls j$ which appears lower in $\tts$ than in $\ttt^\la$, so $\tts$ is not $\la$-row-dominated on $1,\dots,j$.
\end{pf}

\begin{cory}\label{domdom}
Suppose $\la,\mu\in\mptn ln$, and $\tts,\ttt\in\std\mu$.
\begin{enumerate}
\item\label{domdom1}
If $\tts$ is $\la$-dominated on $1,\dots,j$ and $\tts\dom\ttt$, then $\ttt$ is $\la$-dominated on $1,\dots,j$. In particular, if $\tts\in\cstd\la\mu$ and $\tts\dom\ttt$, then $\ttt\in\cstd\la\mu$.
\item\label{domdom2}
If $\tts$ is $\la$-row-dominated on $1,\dots,j$ and $\tts\domby\ttt$, then $\ttt$ is $\la$-row-dominated on $1,\dots,j$. In particular, if $\tts\in\rstd\la\mu$ and $\tts\domby\ttt$, then $\ttt\in\rstd\la\mu$.
\end{enumerate}
\end{cory}

\begin{lemma}\label{wdombru}
Suppose $\la,\mu\in\mptn ln$, and $\ttt,\ttu\in\std\mu$ with $\ttu\domby \ttt$. If $\ttt$ is weakly $\la$-dominated on $1,\dots,j$, then so is $\ttu$.
\end{lemma}
\begin{pf}
The proof follows almost identically to that of \cref{domdom}(\ref{domdom1}).
\end{pf}

\subsection{Dominated homomorphisms}\label{domhsec}

Given $\la,\mu\in\mptn ln$, we want to consider the space of $\hhh n$-homomorphisms $\phi:\spe\la\to\spe\mu$ with the property that $\phi(z_\la)$ lies in the $\bbf$-span of $\lset{v_\tts}{\tts\in\cstd\la\mu}$. But we need to show that this notion is well-defined.

\begin{propn}\label{dominatedbasis}
Suppose $\la,\mu\in\mptn ln$. Then the subspace $\lspan{v_\tts}{\tts\in\cstd\la\mu}_\bbf$ of $\spe\mu$ is independent of the choice of standard basis elements $v_\tts$.
\end{propn}

\begin{pf}
Let $V$ denote the space $\lspan{v_\tts}{\tts\in\cstd\la\mu}_\bbf$, and take $\ttt\in\cstd\la\mu$. Let $s_{j_1}\dots s_{j_r}$ be a new reduced expression for $w_\ttt$, and let $v'_\ttt=\psi_{j_1}\dots\psi_{j_r}z_\mu$ (where $\psi_1,\dots,\psi_{n-1}$ are taken to lie in $\hhh{\cont(\la)}$). Let $V'$ be the space obtained from $V$ by replacing $v_\ttt$ with $v'_\ttt$ in the spanning set $\lset{v_\tts}{\tts\in\cstd\la\mu}$; it suffices to show that $V=V'$. By \cref{diffredex},
\[
v'_\ttt= v_\ttt + \sum_{\ttu\domsby\ttt} a_\ttu v_\ttu \quad \text{for some $a_\ttu\in\bbf$.}
\]
By \cref{domdom}(\ref{domdom1}), each $v_\ttu$ with $\ttu\domby\ttt$ lies in $V$, and so $v'_\ttt\in V$. Hence $V'\subseteq V$; but since the elements $v_\tts$ are linearly independent, $\dim_\bbf V=\dim_\bbf V'=|\cstd\la\mu|$. So $V'=V$.
\end{pf}

In view of \cref{dominatedbasis} and an analogue for row-dominated tableaux, the following definition makes sense.

\begin{defn}
Suppose $\la,\mu\in\mptn ln$. If $\phi\in\shom\la\mu$, we say that $\phi$ is \emph{(column-)dominated} if $\phi(z_\la)\in\lspan{v_\tts}{\tts\in\cstd\la\mu}_\bbf$. We write $\chom\la\mu$ for the space of dominated homomorphisms from $\spe\la$ to $\spe\mu$.

Similarly, if $\chi\in\rshom\la\mu$, we say that $\chi$ is \emph{row-dominated} if $\chi(z^\la)\in\sspan{v^\tts}{\tts\in\rstd\la\mu}_\bbf$, and we write $\rhom\la\mu$ for the space of row-dominated homomorphisms from $\rspe\la$ to $\rspe\mu$.
\end{defn}

\begin{propn}
$\chom\la\mu$ and $\rhom\la\mu$ are graded subspaces of $\shom\la\mu$ and $\rshom\la\mu$ respectively. That is, $\chom\la\mu$ and $\rhom\la\mu$ are spanned by homogeneous homomorphisms.
\end{propn}

\begin{pf}
The proof proceeds almost identically to the proof of the fact that $\shom\la\mu$ is graded, using the fact that $\lspan{v_\tts}{\tts\in\cstd\la\mu}_\bbf$ is a graded subspace of $\spe\mu$ (with a corresponding statement for $\rspe\mu$).
\end{pf}

The rest of this section is devoted to showing that in certain cases every Specht homomorphism is dominated. Specifically, we shall prove the following.

\begin{thm}\label{domhom}
Suppose $e\neq2$ and that $\kappa_1,\dots,\kappa_l$ are distinct. Then $\hom_{\hhh n}(\spe\la,\spe\mu)=\chom\la\mu$.
\end{thm}

\begin{rmk}
The hypotheses that $e\neq2$ and that $\kappa_1,\dots,\kappa_l$ are equivalent to the condition that $\hhh n$ has exactly $2l$ isomorphism classes of one-dimensional modules. These hypotheses also appear in Rouquier's work \cite[Theorems 6.6, 6.8, 6.13]{rouq1} on $1$-faithful quasi-hereditary covers of cyclotomic Hecke algebras. The following small examples show that these hypotheses are essential in \cref{domhom}; in fact, they show that Specht modules labelled by different multipartitions can be isomorphic without these assumptions.
\begin{enumerate}
\item
Take $e=2$, $\kappa=(0)$, $\la=((1^2))$ and $\mu=((2))$. Then there is a non-zero homomorphism $\spe\la\to\spe\mu$ defined by $z_\la\mapsto z_\mu$, though the tableau $\ttt_\mu=\young(12)$ is not $\la$-dominated. So $\shom\la\mu\neq\{0\}=\chom\la\mu$.
\item
For any $e$, take $\kappa=(0,0)$, $\la=(\varnothing,(1))$ and $\mu=((1),\varnothing)$. Then $z_\la\mapsto z_\mu$ again defines a non-zero homomorphism $\spe\la\to\spe\mu$, though $\ttt_\mu$ is not $\la$-dominated.
\end{enumerate}
\end{rmk}

The proof of \cref{domhom} requires several preliminary results. We fix $\la,\mu\in\mptn ln$ and an \emc $\kappa$ of level $l$ throughout. If $\cont(\la)\neq\cont(\mu)$, then by \cref{samedef} $\shom\la\mu=0$, so that \cref{domhom} is trivially true. So we assume that $\cont(\la)=\cont(\mu)$. In the results below, $\psi_1,\dots,\psi_{n-1}$ are elements of $\hhh{\cont(\la)}$.

\begin{lemma}\label{psjdom}
Suppose $j\in\{2,\dots,n\}$ with $j-1\downarrow_{\ttt_\la}j$, and $\ttt\in\std\mu$ is $\la$-dominated on $1,\dots,j$. Then $\psi_{j-1}v_\ttt$ is a linear combination of basis elements $v_\ttu$ for standard tableaux $\ttu$ which are $\la$-dominated on $1,\dots,j$.
\end{lemma}

\begin{pf}
If $j-1\to_\ttt j$ or $j-1\downarrow_\ttt j$ or $j-1\Swarrow_\ttt j$, then the result follows from \cref{domdom}(\ref{domdom1}) together with either \cref{samecol} or \cref{unswitch}. The remaining possibility is that $j-1\Nearrow_\ttt j$. But now if we let $\tts$ be the standard tableau $s_{j-1}\ttt$, then by \cref{diffredex} $\psi_{j-1}v_\ttt = v_\tts + \sum_{\ttu\domsby\tts} b_\ttu v_\ttu$ for some $b_\ttu\in\bbf$. Clearly since $\ttt$ is $\la$-dominated on $1,\dots,j$ and $j-1,j$ lie in the same column of $\ttt_\la$, $\tts$ is also $\la$-dominated on $1,\dots,j$. \cref{domdom}(\ref{domdom1}) completes the proof.
\end{pf}

\begin{propn}\label{jdom}
Suppose $e\neq2$, and that $\phi:\spe\la\to\spe\mu$ is a homomorphism, and write
\[
\phi(z_\la)=\sum_{\ttt\in\std\mu}a_\ttt v_\ttt \quad \text{for some } a_\ttt\in\bbf.
\]
Suppose $j\in\{2,\dots,n\}$ with $j-1\downarrow_{\ttt_\la}j$, and that each $\ttt$ for which $a_\ttt\neq0$ is $\la$-dominated on $1,\dots,j-1$. Then each $\ttt$ for which $a_\ttt\neq0$ is $\la$-dominated on $1,\dots,j$.
\end{propn}

\begin{pf}
The fact that $j-1\downarrow_{\ttt_\la}j$ means that $\psi_{j-1} z_\la=0$, so we must have $\sum_{\ttt\in\std\mu}a_\ttt\psi_{j-1}v_\ttt=0$. Assuming the \lcnamecref{jdom} is false, there is at least one $\ttt$ which is not $\la$-dominated on $1,\dots,j$ such that $a_\ttt\neq0$; choose such a $\ttt$ which is $\dom$-maximal. Since $\ttt$ is $\la$-dominated on $1,\dots,j-1$, the entry $j$ lies in a column strictly to the right of $j-1$ in $\ttt$. We claim that we cannot have $j-1\to_\ttt j$. If this is the case, then the residue sequence $i(\ttt)$ satisfies $i(\ttt)_j=i(\ttt)_{j-1}+1$. However, since $f$ is a homomorphism and $v_\ttt$ appears with non-zero coefficient in $\phi(z_\la)$, we must have $i(\ttt)=i_\la$, and the fact that $j-1\downarrow_{\ttt_\la}j$ means that $(i_\la)_j=(i_\la)_{j-1}-1$. Since $e\neq2$, this is a contradiction.

Hence $j-1\Nearrow_\ttt j$, so the tableau $\tts:=s_{j-1}\ttt$ is standard, and if we write $\psi_{j-1}v_\ttt$ as a linear combination of standard basis elements, then $v_\tts$ occurs with coefficient $1$. We claim that $v_\tts$ does not occur in any other $\psi_{j-1}v_{\ttt'}$ when $a_{\ttt'}\neq0$: if $\ttt'$ is not $\la$-dominated on $1,\dots,j$, then (defining $\tts'$ analogously to $\tts$) we have $\psi_{j-1}v_{\ttt'}=v_{\tts'} + \sum_{\ttu\domsby\tts'} c_\ttu v_\ttu $ for some $c_\ttu\in\bbf$; but the fact that $\ttt\ndomby\ttt'$ (by our choice of $\ttt$ being $\dom$-maximal) means that $\tts\ndomby\tts'$, so $v_\tts$ cannot occur. On the other hand, if $\ttt'$ is $\la$-dominated on $1,\dots,j$, then the result follows from \cref{psjdom}, since $\tts$ is not $\la$-dominated on $1,\dots,j$.

So $v_\tts$ occurs with non-zero coefficient in $\sum_{\ttt\in\std\mu}a_\ttt\psi_{j-1}v_\ttt$, a contradiction.
\end{pf}

We now turn our attention to the case where $j$ is in the top row of its component in $\ttt_\la$.

\begin{lemma}\label{psjwdom}
Suppose $1\ls a\ls j\ls n$, and that $j-1 \nearrow_{\ttt_\la}j$ and $a \downarrow_{\ttt_\la} a+i$ for all $i=1,\dots,j-a-1$. If $\ttt\in\std\mu$ is weakly $\la$-dominated on $1,\dots,j$ then $\psi_a\psi_{a+1}\dots\psi_{j-1}v_\ttt$ is a linear combination of basis elements $v_\ttu$ for standard tableaux $\ttu$ which are weakly $\la$-dominated on $1,\dots,j$.
\end{lemma}

\begin{pf}
We argue by induction on $l(s_a s_{a+1}\dots s_{j-1})=j-a$. If $j-a=0$, the result is trivial. So suppose $a<j$, and assume by induction that $\psi_{a+1}\dots\psi_{j-1}v_\ttt$ is a linear combination of basis elements $v_\ttu$ which are weakly $\la$-dominated on $1,\dots,j$. We want to show that for each $v_\ttu$, $\psi_av_\ttu$ is a linear combination of basis elements $v_{\ttu'}$ for standard tableaux $\ttu'$ which are weakly $\la$-dominated on $1,\dots,j$.

If $a\to_\ttu a+1$ or $a\downarrow_\ttu a+1$ or $a\Swarrow_\ttu a+1$, then the result follows from \cref{wdombru} together with either \cref{samecol} or \cref{unswitch}. The remaining possibility is that $a\Nearrow_\ttu a+1$. Let $\tts$ be the standard tableau $s_a\ttu$. Then by \cref{diffredex}, $\psi_a v_\ttu = v_\tts + \sum_{\ttu'\domsby\tts}a_{\ttu'}v_{\ttu'}$ for some $a_{\ttu'}\in\bbf$.

Recalling that $\ttu$ is weakly $\la$-dominated on $1,\dots,j$ and that $a,a+1$ are in the same column (and therefore the same component) of $\ttt_\la$, $\tts$ is weakly $\la$-dominated on $1,\dots,j$ and \cref{wdombru} completes the proof.
\end{pf}

\begin{propn}\label{wjdom}
Suppose $\phi:\spe\la\to\spe\mu$ is a homomorphism with
\[
\phi(z_\la)=\sum_{\ttt\in\std\mu}a_\ttt v_\ttt \quad \text{for some } a_\ttt\in\bbf.
\]
Suppose $j\in\{2,\dots,n\}$ with either $j-1\nearrow_{\ttt_\la}j$ or $j-1\to_{\ttt_\la}j$, and that each $\ttt$ for which $a_\ttt\neq0$ is $\la$-dominated on $1,\dots,j-1$. Then each $\ttt$ for which $a_\ttt\neq0$ is $\la$-dominated on $1,\dots,j$.
\end{propn}

\begin{pf}
The proof follows the same lines as \cref{jdom}. The condition that $j-1\nearrow_{\ttt_\la}j$ or $j-1\to_{\ttt_\la}j$ means that $\spe\la$ satisfies a Garnir relation $\psi_a \psi_{a+1} \dots \psi_{j-1} z_\la =0$, where $a$ is the entry immediately to the left of $j$ in $\ttt_\la$; since $f$ is a homomorphism, we therefore have $\sum_{\ttt\in\std\mu}a_\ttt\psi_a \dots \psi_{j-1}v_\ttt=0$. Assuming the result is false, there is at least one $\ttt$ which is not $\la$-dominated on $1,\dots,j$ such that $a_\ttt\neq0$; choose such a $\ttt$ which is $\dom$-maximal. Since $\ttt$ is $\la$-dominated on $1,\dots,j-1$ but not $1,\dots,j$, we have $j-1\Nearrow_\ttt j$. In fact $j-1$ and $j$ are in different components of $\ttt$: if not, what is the entry immediately to the left of $j$ in $\ttt$? It must be some $k<j$, since $\ttt$ is standard, but by assumption $k$ is strictly left of $j$ in $\ttt_\la$ and hasn't moved to the right in $\ttt$.

Let $\tts$ denote the standard tableau $s_a s_{a+1} \dots s_{j-1}\ttt$. Then $l(w_\tts)=l(w_\ttt)+j-a$, so that when we write $\psi_a \psi_{a+1} \dots \psi_{j-1}v_\ttt$ as a linear combination of standard basis elements, $v_\tts$ occurs with coefficient $1$. We claim that $v_\tts$ does not occur with non-zero coefficient in $\psi_a \psi_{a+1}\dots\psi_{j-1}v_{\ttt'}$ for any other $\ttt'$ with $a_{\ttt'}\neq0$: if $\ttt'$ is not $\la$-dominated on $1,\dots,j$, then (defining $\tts'$ analogously to $\tts$) we have $\psi_a \psi_{a+1} \dots \psi_{j-1}v_{\ttt'}=v_{\tts'} + \sum_{\ttu\domsby\tts'} b_\ttu v_\ttu$ for some $b_\ttu\in\bbf$; but the fact that $\ttt\ndomby\ttt'$ (by our choice of $\ttt$) means that $\tts\ndomby\tts'$, so $v_\tts$ cannot occur. On the other hand, if $\ttt'$ is $\la$-dominated on $1,\dots,j$, then the result follows from \cref{psjwdom}, since $\tts$ is not weakly $\la$-dominated on $1,\dots,j$ as $j-1$ and $j$ are in different components of $\ttt$.

So $v_\tts$ occurs with non-zero coefficient in $\sum_{\ttt\in\std\mu}a_\ttt\psi_a\psi_{a+1}\dots\psi_{j-1}v_\ttt$, a contradiction.
\end{pf}

The last thing we need for the proof of \cref{domhom} is the following.

\begin{lemma}\label{topleft}
Suppose $\kappa_1,\dots,\kappa_l$ are distinct, and that $\ttt\in\std\mu$ satisfies $i(\ttt)=i_\la$. If $\ttt$ is $\la$-dominated on $1,\dots,j-1$ and $j$ appears in the $(1,1)$-position of its component in $\ttt_\la$, then $\ttt$ is $\la$-dominated on $1,\dots,j$.
\end{lemma}

\begin{pf}
Suppose not; then $j$ appears in $\ttt$ strictly to the right of where it appears in $\ttt_\la$. This means that $j$ must appear in the $(1,1)$-node of some component of $\ttt$, since otherwise there would be a smaller entry immediately above or to the left of $j$, contradicting the assumption that $\ttt$ is $\la$-dominated on $1,\dots,j-1$.

So there are $1\ls r<s\ls l$ such that $\ttt_\la(1,1,s)=j=\ttt(1,1,r)$. Hence $\kappa_s=(i_\la)_j=i(\ttt)_j=\kappa_r$, contrary to assumption.
\end{pf}

\begin{pf}[Proof of \cref{domhom}]
Suppose $\phi:\spe\la\to\spe\mu$ is a homomorphism, and write
\[
\phi(z_\la)=\sum_{\ttt\in\std\mu}a_\ttt v_\ttt \quad \text{for some }a_\ttt\in\bbf.
\]
We must show that every $\ttt$ for which $a_\ttt\neq0$ is $\la$-dominated. In fact we show by induction on $j$ that every such $\ttt$ is $\la$-dominated on $1,\dots,j$, with the case $j=0$ being vacuous. So suppose $j\gs1$, and assume by induction that $\ttt$ is $\la$-dominated on $1,\dots,j-1$. Note that since $\phi$ is a homomorphism, we have $i(\ttt)=i_\la$.

If $j=1$ or $j$ lies in an earlier component of $\ttt_\la$ than $j-1$, then $j$ lies in the $(1,1)$-node of its component in $\ttt_\la$. So by \cref{topleft} $\ttt$ is $\la$-dominated on $1,\dots,j$. The remaining possibilities are that $j>1$ and that one of
\[
j-1\downarrow_{\ttt_\la}j,\qquad j-1\to_{\ttt_\la}j,\qquad j-1\nearrow_{\ttt_\la}j
\]
occurs; these cases are dealt with in \cref{jdom,wjdom}.
\end{pf}

We immediately see the following interesting result.

\begin{cory}
Suppose $e\neq2$ and that $\kappa_1,\dots,\kappa_l$ are distinct. If $\la,\mu\in\mptn ln$ with $\Hom_{\hhh n}(\spe\la,\spe\mu)\neq\{0\}$, then $\la\dom\mu$. Furthermore (since $\cstd\la\la=\{\ttt_\la\}$)  $\hom_{\hhh n}(\spe\la,\spe\la)$ is one-dimensional. In particular, $\spe\la$ is indecomposable.
\end{cory}

\begin{rmk}
Note that if $e=2$ then $\spe\la$ may be decomposable. For example, when $l=1$ and $\nchar(\bbf)\neq3$, the Specht module $\spe{((5,1^2))}$ is decomposable; this was shown in \cite[Example 23.10(iii)]{jbook} in the case $\nchar(\bbf)=2$, and in \cite[Theorem 6.8]{ls} in odd characteristic. Similarly, when $\kappa_i=\kappa_j$ for some $i\neq j$, we can have decomposable Specht modules: take $\kappa=(0,0)$, $e=3$ and $\nchar(\bbf)\neq2$; then $\spe{((3),(3))}$ is decomposable.
\end{rmk}

In exactly the same way, we can prove the corresponding result for row Specht modules.

\begin{thm}\label{domhomrow}
Suppose $e\neq2$ and that $\kappa_1,\dots,\kappa_l$ are distinct, and $\la,\mu\in\mptn ln$. Then $\rhom\la\mu=\rshom\la\mu$. Hence $\rshom\la\mu\neq\{0\}$ only if $\la\domby\mu$, $\rshom\la\la$ is one-dimensional, and $\rspe\la$ is indecomposable.
\end{thm}

\subsection{Duality for dominated homomorphisms}\label{dualsec}

In this section we consider the relationship between row and column Specht modules, as well as the relationship between Specht modules labelled by conjugate multipartitions. These relationships are encapsulated in \cite[Theorems 7.25 and 8.5]{kmr}, from which it follows that a (generalised) column-removal theorem for homomorphisms between Specht modules is equivalent to the corresponding row-removal theorem. The main result of this section, which requires considerable additional work, is that the same is true for dominated homomorphisms.

Following \cite[\S3.2]{kmr}, let $\tau:\hhh\alpha\to\hhh\alpha$ denote the anti-automorphism which fixes all the generators $e(i),y_r,\psi_s$, and define $\tau:\hhh n\to\hhh n$ by combining these maps for all $\alpha$. If $M=\bigoplus_{d\in\bbz}M_d$ is a graded $\hhh n$-module, let $M\dual$ denote the graded module with $M\dual_d=\hom_\bbf(M_{-d},\bbf)$ for each $d$, with $\hhh n$-action given by $(hf)m=f(\tau(h)m)$ for $m\in M$, $f\in M\dual$ and $h\in \hhh n$. Also, for $k\in\bbz$, let $M\lan k\ran$ denote the same module with the grading shifted by $k$, i.e.\ $M\lan k\ran_d=M_{d-k}$. Finally, recall the defect $\df(\la)$ of a multipartition from \cref{resdegsec}.

\begin{thmc}{kmr}{Theorem 7.25}\label{7.25}
Suppose $\la\in\mptn ln$. Then
\[
\rspe\la\cong(\spe\la)\dual\lan\df(\la)\ran\qquad\text{and}\qquad\spe\la\cong(\rspe\la)\dual\lan\df(\la)\ran.
\]
\end{thmc}

Now suppose $\la,\mu\in\mptn ln$. Applying \cref{7.25} to both $\la$ and $\mu$ gives an isomorphism of graded vector spaces
\[
\rshom\mu\la\cong\hom_{\hhh n}(\spe\mu\dual\lan\df(\mu)\ran,\spe\la\dual\lan\df(\la)\ran);
\]
since by \cref{samedef} $\df(\la)=\df(\mu)$ for any $\la$ and $\mu$ with $\rshom\mu\la\neq\{0\}$, this yields an isomorphism of graded vector spaces
\[
\rshom\mu\la\cong\hom_{\hhh n}(\spe\mu\dual,\spe\la\dual).
\]
$\tau$ is a homogeneous map of degree zero, so $\hom_{\hhh n}(\spe\mu\dual,\spe\la\dual)$ is canonically isomorphic as a graded vector space to $\shom\la\mu$, and hence we have an isomorphism of graded vector spaces
\[
\Theta:\shom\la\mu\stackrel\sim\longrightarrow\rshom\mu\la.
\]
Our aim is to prove the following.

\begin{propn}\label{rowcold}
Suppose $\la,\mu\in\mptn ln$, and let $\Theta:\shom\la\mu\to\rshom\mu\la$ be the bijection above. Then $\Theta(\chom\la\mu)=\rhom\mu\la$.
\end{propn}

We shall prove \cref{rowcold} below. First we examine the consequences for row and column removal. In order to be able to compare row and column removal, we combine \cref{rowcold} with a result which relates to an analogue of the sign representation of the symmetric group. Following \cite[\S3.3]{kmr}, let $\sgn:\hhh\alpha\to\hhh\alpha$ denote the automorphism which maps $e(i)\mapsto e(-i)$, $y_r\mapsto-y_r$ and $\psi_s\mapsto-\psi_s$ for all $i,r,s$, and define $\sgn:\hhh n\to\hhh n$ by combining these maps for all $\alpha$. Given a graded $\hhh n$-module $M$, let $M^{\sgn}$ denote the same graded vector space with the action of $\hhh n$ twisted by $\sgn$.

Recall that if $\la$ is a multipartition, then $\la'$ denotes the conjugate multipartition to $\la$, and that if $\tts\in\std\la$, then $\tts'\in\std{\la'}$ denotes the conjugate tableau to $\tts$. Also define the \emph{conjugate \emc} $\kappa':=(-\kappa_l,\dots,-\kappa_1)$. Now the following is immediate from the construction of row and column Specht modules.

\begin{thmc}{kmr}{Theorem 8.5}\label{sign}
Suppose $\la\in\mptn ln$. Then there is an isomorphism $(\rspek\la\kappa)^{\sgn}\cong\spek{\la'}{\kappa'}$ of $\hhh n$-modules, given by $v^\tts\mapsto v_{\tts'}$.
\end{thmc}

\begin{rmk}
\cref{sign} is one place where it is essential that we consider Specht modules as modules for $\hhh n$, rather than its cyclotomic quotients, since the two modules involved are defined relative to different \emcs.
\end{rmk}

Now suppose $\la,\mu\in\mptn ln$. Since $\sgn$ is a homogeneous automorphism of $\hhh n$, we have an equality of graded vector spaces
\begin{align*}
\hom_{\hhh n}((\rspek\mu\kappa)^{\sgn},(\rspek\la\kappa)^{\sgn})&=\hom_{\hhh n}(\rspek\mu\kappa,\rspek\la\kappa)\tag{\textasteriskcentered},
\\
\intertext{Combining this with \cref{sign}, we have an isomorphism of graded vector spaces}
\hom_{\hhh n}(\spek{\mu'}{\kappa'},\spek{\la'}{\kappa'})&\cong\hom_{\hhh n}(\rspek\mu\kappa,\rspek\la\kappa)\tag*{(\textdagger)}.
\\
\intertext{Applying \cref{7.25} yields an isomorphism of graded vector spaces}
\hom_{\hhh n}(\spek{\mu'}{\kappa'},\spek{\la'}{\kappa'})&\cong\hom_{\hhh n}(\spek\la\kappa,\spek\mu\kappa)\tag*{(\textdaggerdbl)}.
\end{align*}
We want to show that the same holds for dominated homomorphisms; this is immediate when $e>2$ and $\kappa_1,\dots,\kappa_l$ are distinct, by \cref{domhom}. In general, we observe that (\textasteriskcentered) remains true with $\hom$ replaced by $\DHom$, and the explicit form of the isomorphism in \cref{sign} shows that (\textdagger) does too, since $\tts\in\cstd{\mu'}{\la'}$ if and only if $\tts'\in\rstd{\mu}{\la}$. Finally, \cref{rowcold} shows that (\textdaggerdbl) remains true for $\DHom$ too. So we have the following theorem.

\begin{thm}\label{homconj}
Suppose $\la,\mu\in\mptn ln$. Then there is an isomorphism of graded vector spaces
\[
\dhom_{\hhh n}(\spek\la\kappa,\spek\mu\kappa)\cong\dhom_{\hhh n}(\spek{\mu'}{\kappa'},\spek{\la'}{\kappa'}).
\]
\end{thm}

It remains to prove \cref{rowcold}; for the remainder of this section, all Specht modules are defined for the \emc $\kappa$.

We begin by recalling how the isomorphism $\rspe\la\cong\spe\la\dual\lan\df(\la)\ran$ in \cref{7.25} is constructed. Given the standard basis $\lset{v_\ttt}{\ttt\in\std\la}$ for $\spe\la$, let $\set{f^\ttt}{\ttt\in\std\la}$ be the dual basis for $\spe\la\dual$; although the elements $f^\ttt$ in general depend on the choice of the elements $v_\ttt$ (i.e.\ on the choice of preferred reduced expressions), it is an easy exercise to show that the element $f^{\ttt^\la}$ does not. The isomorphism $\rspe\la\to\spe\la\dual\lan\df(\la)\ran$ is defined by $z^\la\mapsto f^{\ttt^\la}$.

\begin{lemma}\label{triang}
Suppose $\la\in\mptn ln$, and let $\theta^\la:\rspe\la\to\spe\la\dual\lan\df(\la)\ran$ be the isomorphism constructed above.
\begin{enumerate}
\item\label{triang1}
For any $\tts\in\std\la$ we have $\theta^\la(v^\tts)\in\sspan{f^\ttt}{\ttt\in\std\la,\ \ttt\dom\tts}_\bbf$.
\item\label{triang2}
$\theta^\la$ maps the space $\sspan{v^\tts}{\tts\in\rstd\mu\la}_\bbf$ bijectively to the space $\sspan{f^\tts}{\tts\in\rstd\mu\la}_\bbf$.
\end{enumerate}
\end{lemma}

\begin{pfenum}
\item
For each $\ttt\in\std\la$, write $\tau(\psi^\tts) v_\ttt=\sum_{\ttu\in\std\la}a_{\ttt\ttu}v_\ttu$. Then one can check that the definitions give $\theta^\la(v^\tts)=\sum_{\ttt\in\std\la}a_{\ttt\ttt^\la} f^\ttt$. So it suffices to show that $a_{\ttt\ttt^\la}=0$ when $\ttt\ndom\tts$. Clearly to prove this it is sufficient to show this in the case where $\bbf=\bbc$, and so (as in the proof of \cite[Theorem 7.25]{kmr}) we can invoke the proof of \cite[Proposition 6.19]{hm}; here $\theta^\la$ is given in the form $x\mapsto\{x,-\}$, for a bilinear form $\{\ ,\ \}:\rspe\la\times\spe\la\lan\df(\la)\ran\to\bbc$ satisfying $\{v^\tts,v_\ttt\}=0$ unless $\ttt\dom\tts$, which is exactly what we want.
\item
From (\ref{triang1}) and \cref{domdom}(\ref{domdom2}) we have $\theta^\la(v^\tts)\in\sspan{f^\ttt}{\ttt\in\rstd\mu\la}_\bbf$ whenever $\tts\in\rstd\mu\la$, so $\theta^\la\left(\sspan{v^\tts}{\tts\in\rstd\mu\la}\right)\subseteq\sspan{f^\tts}{\tts\in\rstd\mu\la}$. But $\theta^\la$ is an isomorphism of vector spaces and
\[
\dim_\bbf\sspan{v^\tts}{\tts\in\rstd\mu\la}_\bbf=\left|\rstd\mu\la\right|=\dim_\bbf\sspan{f^\tts}{\tts\in\rstd\mu\la}_\bbf,
\]
so in fact $\theta^\la\left(\sspan{v^\tts}{\tts\in\rstd\mu\la}_\bbf\right)=\sspan{f^\tts}{\tts\in\rstd\mu\la}_\bbf$.
\end{pfenum}

\begin{lemma}\label{rowdomlemma}
Suppose $\la,\mu\in\mptn ln$. Suppose $\tts\in\std\la$ and $\ttu$ is a $\la$-tableau such that $w_\tts\cgs w_\ttu$ and that for every $1\ls i\ls n$ the number $i$ appears in $\ttu$ weakly to the right of where it appears in $\ttt^\mu$. Then $\tts\in\rstd\mu\la$.
\end{lemma}

\begin{pf}
Using \cref{domalt}(\ref{domalt2}) we just need to show that $\shp\tts m\dom\shp{\ttt^\mu} m$ for all $m$. Let $\ttu\cstr$ be the column-strict tableau which is column-equivalent to $\ttu$. Then by \cref{tbruhat} $w_\ttu\cgs w_{\ttu\cstr}$. By \cref{brudomconj}, we have that $\shp{(\ttu\cstr)} m'\dom\shp\tts m'$ for all $m$. Furthermore, the condition that every entry in $\ttu\cstr$ lies weakly to the right of where it lies in $\ttt^\mu$ is equivalent to every entry in $(\ttu\cstr)'$ lying weakly below where it lies in $(\ttt^\mu)'$, so we necessarily have that $\shp{(\ttu\cstr)}m'\domby\shp{\ttt^\mu} m'$ for all $m$. Reapplying \cref{brudomconj}, we have $w_{\ttu\cstr}\cgs w_{\ttt^\mu}$.
\end{pf}

\begin{lemma}\label{dualitylemma}
Suppose $\la,\mu\in\mptn ln$, $\tts\in\std\la\setminus\rstd\mu\la$ and $\ttt\in\cstd\la\mu$. Then when $\psi_\tts v_\ttt$ is expressed in terms of the standard basis $\lset{v_\ttu}{\ttu\in\std\mu}$, the coefficient of $v_{\ttt^\mu}$ is zero.
\end{lemma}

\begin{pf}
Suppose to the contrary that $v_{\ttt^\mu}$ does appear with non-zero coefficient in $\psi_\tts v_\ttt=\psi_\tts\psi_\ttt z_\mu$. Let $s_{i_1}\dots s_{i_a}$ and $s_{j_1}\dots s_{j_b}$ be the preferred reduced expressions for $w_\tts$ and $w_\ttt$ respectively. Then by \cref{subexp} there is a reduced expression for $w_{\ttt^\mu}$ occurring as a subexpression of $s_{i_1}\dots s_{i_a}s_{j_1}\dots s_{j_b}$. If we separate this reduced expression into two parts, which occur as subexpressions of $s_{i_1}\dots s_{i_a}$ and $s_{j_1}\dots s_{j_b}$ respectively, and let $w,x$ denote the corresponding elements of $\sss n$, then we have
\[
w\cls w_\tts,\ x\cls w_\ttt,\ wx=w_{\ttt^\mu},\ l(w)+l(x)=l(w_{\ttt^\mu}).
\]
Putting $\ttv=x\ttt_\mu$, we have $\ttv\in\std\mu$ by \cref{leftstd}, and in fact $\ttv\in\cstd\la\mu$ (using \cref{domdom}(\ref{domdom1}), because $w_\ttv\cls w_\ttt$ and $\ttt\in\cstd\la\mu$). If we let $\ttu=w\ttt_\la$ then, as functions $[\mu]\to[\la]$,
\[
\ttu^{-1}\ttt^\mu=\ttt_\la^{-1}x\ttt_\mu=\ttt_\la^{-1}\ttv.
\]
The fact that $\ttv$ is $\la$-dominated can be expressed as saying that the map $\ttt_\la^{-1}\ttv:[\mu]\to[\la]$ maps any node of $\mu$ to a node weakly to the right. So each entry of $\ttu$ appears weakly to the right of where it appears in $\ttt^\mu$, i.e.\ $\ttu$ satisfies the hypotheses of \cref{rowdomlemma}. Hence by \cref{rowdomlemma} $\tts\in\rstd\mu\la$, contrary to hypothesis.
\end{pf}

\begin{pf}[Proof of \cref{rowcold}]
We shall prove that $\Theta(\chom\la\mu)\subseteq\rhom\mu\la$; the same argument with $\la$ and $\mu$ interchanged and with row and column Specht modules interchanged proves the opposite containment.

Suppose $\phi\in\chom\la\mu$, and write $\phi(z_\la)=\sum_{\ttt\in\cstd\la\mu}a_\ttt v_\ttt$ for some $a_\ttt\in\bbf$. Let $\phi\dual:\spe\mu\dual\to\spe\la\dual$ denote the dual map. We want to show that the homomorphism $\Theta(\phi)$ which corresponds to $\phi\dual$ via \cref{7.25} is row-dominated, i.e.{} $\Theta(\phi)(z^\mu)\in\sspan{v^\tts}{\tts\in\rstd\mu\la}_\bbf$. By the construction of the isomorphism $\rspe\mu\to\spe\mu\dual$ and by \cref{triang}, this is the same as saying that $\phi\dual(f^{\ttt^\mu})\in\sspan{f^\tts}{\tts\in\rstd\mu\la}_\bbf$; in other words, $\phi\dual(f^{\ttt^\mu})(v_\tts)=0$ when $\tts\in\std\la\setminus\rstd\mu\la$.

$\phi\dual$ is given by $f\mapsto f\circ\phi$. In particular, $\phi\dual(f^{\ttt^\mu})=f^{\ttt^\mu}\circ\phi$, which maps $v_\tts$ to the coefficient of $v_{\ttt^\mu}$ in $\phi(v_\tts)=\sum_{\ttt\in\cstd\la\mu}a_\ttt\psi_\tts v_\ttt$. By \cref{dualitylemma} this coefficient is zero when $\tts\notin\rstd\mu\la$, and the result follows.
\end{pf}

\section{Column removal for homomorphisms}\label{crhsec}

Now we come to the main results of the paper, which give row and column removal theorems for dominated homomorphisms between Specht modules.

\subsection{Generalised column removal for multipartitions}\label{gcrnotation}

\begin{defn}
Suppose $\la=(\la^{(1)},\dots,\la^{(l)})\in\mptn ln$. For any $1\leq m\leq l$ and any $c\gs0$, define $\la^{(m),c}\lf$ to be the partition consisting of all nodes in the first $c$ columns of $\la^{(m)}$, and $\la^{(m),c}\rg$ the partition consisting of all nodes after the first $c$ columns of $\la^{(m)}$. That is,
\[
(\la\lf^{(m),c})_i=\min\left\{\la^{(m)}_i,c\right\},\qquad(\la\rg^{(m),c})_i=\max\left\{\la^{(m)}_i-c,0\right\}\tag*{for all $i\gs1$.}
\]
Now define
\begin{align*}
\la\rg&=\la\rg(c,m)=(\la^{(1)},\dots,\la^{(m-1)},\la^{(m),c}\rg),\\
\la\lf&=\la\lf(c,m)=(\la^{(m),c}\lf,\la^{(m+1)},\dots,\la^{(l)}).
\end{align*}
Here is an enlightening pictorial representation of this construction, with $l=3$, $m=2$ and $c=3$.
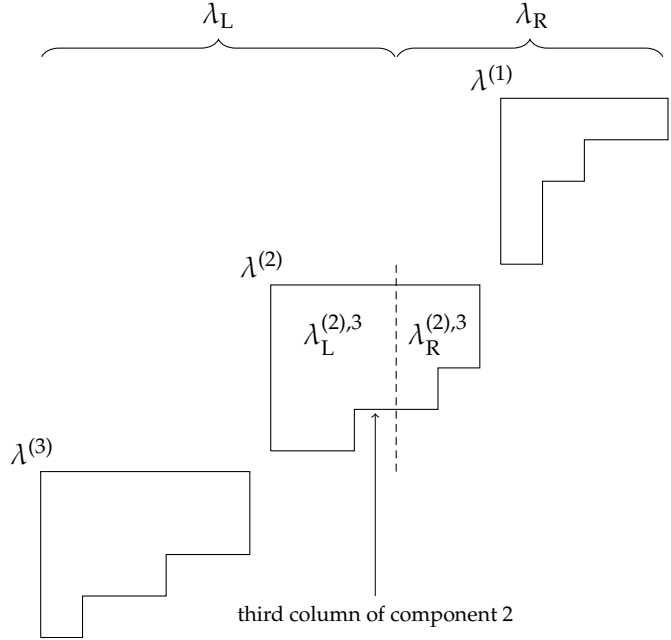
\[
\begin{tikzpicture}[scale=.55,decoration={brace,amplitude=.6em}]
\draw (11,9) -- (11,13) -- (15,13) -- (15,12) -- (13,12) -- (13,11) -- (12,11) -- (12,9) -- (11,9);
\node at (10.8,13.5) {$\la^{(1)}$};

\draw (5.5,4.5) -- (5.5,8.5) -- (10.5,8.5) -- (10.5,6.5) -- (9.5,6.5) -- (9.5,5.5) -- (7.5,5.5) -- (7.5,4.5) -- (5.5,4.5);
\node at (5.3,9) {$\la^{(2)}$};
\node at (7,7.3) {$\la^{(2),3}\lf$};
\node at (9.5,7.3) {$\la^{(2),3}\rg$};

\draw (0,0) -- (0,4) -- (5,4) -- (5,2) -- (3,2) -- (3,1) -- (1,1) -- (1,0) -- (0,0);
\node at (-0.2,4.5) {$\la^{(3)}$};

\draw[densely dashed] (8.5,4) -- (8.5,9);

\draw[->,very thin] (8,1) node[below] {\scriptsize{third column of component 2}} -- (8,5.4);

\draw[decorate](.025,14)--++(8.45,0);\draw(4.25,15)node{$\la\lf$};
\draw[decorate](8.525,14)--++(6.35,0);\draw(11.75,15)node{$\la\rg$};
\end{tikzpicture}
\]
Now we consider tableaux. Suppose $\la\lf,\la\rg$ are as above, and let $n\lf=|\la\lf|$ and $n\rg=|\la\rg|$. Given $\ttt\lf\in\std{\la\lf}$ and $\ttt\rg\in\std{\la\rg}$, define $\lr\ttt\ttt$ to be the $\la$-tableau obtained by filling in the entries $1,\dots,n\lf$ as they appear in $\ttt\lf$, and then filling in the entries $n\lf+1,\dots,n$ as $1,\dots,n\rg$, respectively, appear in $\ttt\rg$. Observe that if $\ttt\in\std\la$ and the integers $1,\dots,n\lf$ all appear in $\ttt$ in column $c$ of component $m$ or further to the left, then $\ttt$ has the form $\lr\ttt\ttt$ for some $\ttt\lf\in\std{\la\lf}$ and $\ttt\rg\in\std{\la\rg}$. We write $\stdlr\la$ for the set of $\ttt\in\std\la$ with this property.
\end{defn}

\begin{eg}
Take $l=3$ and $\la=\left((3),(2^2),(2,1)\right)$. Taking $m=2$ and $c=1$, we get
\[
\la\lf=\left((1^2),(2,1)\right),\qquad\la\rg=\left((3),(1^2)\right).
\]
If we choose
\[
\ttt\lf=\gyoung(::^\hf;1,::^\hf;3,/\hf,;2;4,;5),\qquad\ttt\rg=\gyoung(:^\hf;2;3;5,/\hf,;1,;4),
\]
then we obtain
\[
\lr\ttt\ttt=\gyoung(^5;7;8;\ten,/\hf,::^\hf;1;6,::^\hf;3;9,/\hf,;2;4,;5).
\]
\end{eg}

\subsection{Simple row and column removal}

\begin{thm}[Graded Column Removal]\label{cr}
Suppose $\la,\mu\in\mptn ln$ and $1\ls m\ls l$. Suppose that $\la^{(m+1)}=\dots=\la^{(l)}=\mu^{(m+1)}=\dots=\mu^{(l)}=\varnothing$, and $k:=(\la^{(m)'})_1=(\mu^{(m)'})_1$. Let $\la\rg=\la\rg(1,m)$, $\mu\rg=\mu\rg(1,m)$ and $\kappa\rg=(\kappa_1,\dots,\kappa_{m-1},\kappa_m+1)$. Then
\[
\dhom_{\hhh n}(\spek\la\kappa,\spek\mu\kappa)\cong\DHom_{\hhh{n-k}}(\spek{\la\rg}{\kappa\rg},\spek{\mu\rg}{\kappa\rg})
\]
as graded vector spaces over $\bbf$.
\end{thm}

\begin{rmk}
Recalling \cref{domhom}, this result in fact implies that $\hom_{\hhh n}(\spek\la\kappa,\spek\mu\kappa)\cong\hom_{\hhh{n-k}}(\spek{\la\rg}{\kappa\rg},\spek{\mu\rg}{\kappa\rg})$ when $e\neq2$ and $\kappa_1,\dots,\kappa_l$ are distinct.
\end{rmk}

\begin{pf}
We construct the isomorphism explicitly in the KLR setting. First note that we may assume $\la\dom\mu$, since otherwise $\cstd{\la\rg}{\mu\rg}=\cstd\la\mu=\emptyset$ and the result is immediate. We also observe that $\cont(\la)=\cont(\mu)$ if and only if $\cont(\la\rg)=\cont(\mu\rg)$; if these conditions do not hold then the result is trivial since both homomorphism spaces are zero, so we assume $\cont(\la)=\cont(\mu)$, and set $\alpha:=\cont(\la)$, $\beta:=\cont(\la\rg)$.

For this proof we make an assumption about the choice of preferred reduced expressions defining the standard bases for $\spek{\mu\rg}{\kappa\rg}$ and $\spek\mu\kappa$. Given $\ttt\in\cstd{\la\rg}{\mu\rg}$, we define $\ttt^+:=\ttt_{\mu\lf}\cj\ttt$, where
\[
\mu\lf=\mu\lf(1,m)=\left((1^k),\varnothing,\dots,\varnothing\right)\in\mptn{l-m+1}k.
\]
In other words, $\ttt^+$ is obtained from $\ttt$ by increasing each entry by $k$, adding the column $\gyoung(;1,|\sq\vdts,;k)$ at the left of component $m$, and then adding $l-m$ empty components at the end. Now recall the maps (both denoted $\shift k$) from $\sss{n-k}$ to $\sss n$ and from $\hhh\beta$ to $\hhh\alpha$. Observe that for $\ttt\in\cstd{\la\rg}{\mu\rg}$ we have $w_{\ttt^+}=\shift k(w_\ttt)$. By choosing compatible reduced expressions for $w_{\ttt^+}$ and $w_\ttt$, we may assume that $\psi_{\ttt^+}=\shift k(\psi_\ttt)$ as well.

Now let $c=(\la^{(m)}\rg)'_1$. Then the entries $1,\dots,c$ all appear in the first column of component $m$ in $\ttt_{\la\rg}$, and hence if $\ttt\in\cstd{\la\rg}{\mu\rg}$ these entries all appear in the first column of component $m$ of $\ttt$. In particular, $w_\ttt$ fixes $1,\dots,c$, so $\psi_\ttt$ only involves terms $\psi_j$ for $j>c$; hence $\psi_{\ttt^+}$ only involves terms $\psi_j$ for $j>k+c$.

Now suppose $\varphi\rg\in\DHom_{\hhh{n-k}}(\spek{\la\rg}{\kappa\rg},\spek{\mu\rg}{\kappa\rg})$. Then
\begin{align*}
\varphi\rg(z_{\la\rg})&=\sum_{\ttt\in\cstd{\la\rg}{\mu\rg}}a_\ttt v_\ttt\quad\text{for some $a_\ttt\in\bbf$.}
\\
\intertext{We define $\phi:\spek\la\kappa\to\spek\mu\kappa$ by}
\varphi(z_{\la})&=\sum_{\ttt\in\cstd{\la\rg}{\mu\rg}}a_\ttt v_{\ttt^+}.
\end{align*}
We must verify that this does indeed define a homomorphism, i.e.\ that $h\phi(z_\la)=0$ whenever $h\in\ann(z_\la)$. (Here and henceforth we write $\ann(z_\la)$ for the annihilator of $z_\la$.) Firstly, note that if $\ttt\in\cstd{\la\rg}{\mu\rg}$ with $a_\ttt\neq0$, then $\ttt$ has residue sequence $i_{\la\rg}$; this implies that $\ttt^+$ has residue sequence $i_\la$, so that $e(i_\la)\phi(z_\la)=\phi(z_\la)$, as required. For the other relations, observe from the defining relations for the column Specht module that $\shift k(\ann(z_{\la\rg}))\subseteq\ann(z_\la)$ (and similarly for $\mu\rg$ and $\mu$). Now for $k<j\ls n$ we have $y_{j-k}\in\ann(z_{\la\rg})$, so (since $\phi\rg$ is a homomorphism) $y_{j-k}\sum_\ttt a_\ttt\psi_\ttt\in\ann(z_{\mu\rg})$. Hence
\[
\ann(z_\mu)\ni\shift k\left(y_{j-k}\sum_\ttt a_\ttt\psi_\ttt\right)=y_j\sum_\ttt a_\ttt\psi_{\ttt^+},
\]
so that $y_j\phi(z_\la)=0$. A similar statement applies to $\psi_j$ whenever $k<j<n$ with $j\downarrow_{\ttt_\la}j+1$, and to any Garnir element $\garn A$ where $A$ does not lie in the first column of component $m$.

It remains to check the generators of $\ann(z_\la)$ which do not lie in $\shift k(\ann(z_{\la\rg}))$, i.e.\ the elements $y_1,\dots,y_k$, $\psi_1,\dots,\psi_{k-1}$ and $\garn A$ for $A$ of the form $(j,1,m)$ with $1\ls j\ls c$. Let $h$ denote any of these elements, and observe that since each $\psi_{\ttt^+}$ is a product of terms $\psi_i$ with $i>k+c$, $h$ commutes with $\psi_{\ttt^+}$ (note that if $h=\garn{(j,1,m)}$, then $h$ only involves terms $\psi_i$ for $i<k+c$). Hence
\[
h\phi(z_\la)=h\sum_\ttt a_\ttt\psi_{\ttt^+}z_\mu=\sum_\ttt a_\ttt\psi_{\ttt^+}hz_\mu=0,
\]
since $h\in\ann(z_\mu)$.

So $\ann(z_\la)\phi(z_\la)=0$, and $\phi$ is a well-defined homomorphism. So we have a map $\Phi:\DHom_{\hhh{n-k}}(\spek{\la\rg}{\kappa\rg},\spek{\mu\rg}{\kappa\rg})\rightarrow\DHom_{\hhh n}(\spek\la\kappa,\spek\mu\kappa)$ given by $\phi\rg\mapsto\phi$, and $\Phi$ is obviously linear. To show that $\Phi$ is bijective, we construct its inverse. Any $\tts\in\cstd\la\mu$ must have entries $1,\dots,k$ in order down the first column of its $m$th component; that is, $\tts=\ttt^+$ for some $\ttt\in\cstd{\la\rg}{\mu\rg}$. So given $\theta\in\DHom_{\hhh n}(\spek\la\kappa,\spek\mu\kappa)$, we can write 
\begin{align*}
\theta(z_\la)&=\sum_{\ttt\in\cstd{\la\rg}{\mu\rg}}a_\ttt v_{\ttt^+}\quad\text{for some $a_\ttt\in\bbf$.}
\\
\intertext{Applying (a simpler version of) the above argument in reverse, we see that we have a homomorphism $\theta\rg:\spe{\la\rg}\to\spe{\mu\rg}$ given by}
\theta\rg(z_{\la\rg})&=\sum_{\ttt\in\cstd{\la\rg}{\mu\rg}}a_\ttt v_\ttt.
\end{align*}
So we get a linear map $\DHom_{\hhh n}(\spek\la\kappa,\spek\mu\kappa)\to\DHom_{\hhh{n-k}}(\spek{\la\rg}{\kappa\rg},\spek{\mu\rg}{\kappa\rg})$ which is a two-sided inverse to $\Phi$, and hence $\Phi$ is a bijection.

Finally, to show that we have an isomorphism of graded vector spaces, we show that $\Phi$ is homogeneous of degree $0$. That is, if $0\neq\phi\rg\in\DHom_{\hhh{n-k}}(\spek{\la\rg}{\kappa\rg},\spek{\mu\rg}{\kappa\rg})$ is homogeneous, then $\phi$ is also homogeneous with $\deg(\varphi)=\deg(\varphi\rg)$. To see this, we write
\begin{align*}
\varphi\rg(z_{\la\rg})&=\sum_{\ttt\in\cstd{\la\rg}{\mu\rg}}a_\ttt v_\ttt \quad \text{for some $a_\ttt\in\bbf$.}\\
\intertext{Then}
\varphi(z_{\la})&=\sum_{\ttt\in\cstd{\la\rg}{\mu\rg}}a_\ttt v_{\ttt^+},
\end{align*}
and for each $\ttt$ with $a_\ttt\neq0$ we have
\[
\codeg^\kappa(\ttt^+)-\codeg^\kappa(\ttt_\la)=\codeg^{\kappa_R}(\ttt)-\codeg^{\kappa_R}(\ttt_{\la\rg})=\deg(\varphi\rg).
\]
Hence $\phi$ is homogeneous of degree $\deg(\phi\rg)$.
\end{pf}

Now we make corresponding definitions for row removal.

\begin{defn}
Suppose $\la\in\mptn ln$. For any $1\leq m\leq l$ and any $r\gs0$, define
\[
\la^{(m),r}\tp=(\la^{(m)}_1,\dots,\la^{(m)}_r,0,0,\dots),\qquad\la^{(m),r}\bt=(\la^{(m)}_{r+1},\la^{(m)}_{r+2},\dots).
\]
Now let
\begin{align*}
\la\tp&=\la\tp(r,m)=(\la^{(1)},\dots,\la^{(m-1)},\la^{(m),r}\tp),\\
\la\bt&=\la\bt(r,m)=(\la^{(m),r}\bt,\la^{(m+1)},\dots,\la^{(l)}),
\end{align*}
and set $n\tp=|\la\tp|$ and $n\bt=|\la\bt|$.
\end{defn}

\begin{cory}[Graded Row Removal]\label{rr}
Suppose $\la,\mu\in\mptn ln$ and $1\ls m\ls l$. Suppose that $\la^{(1)}=\dots=\la^{(m-1)}=\mu^{(1)}=\dots=\mu^{(m-1)}=\varnothing$, and $k:=\la^{(m)}_1=\mu^{(m)}_1$. Let $\la\bt=\la\bt(1,m)$, $\mu\bt=\mu\bt(1,m)$ and $\kappa\bt=(\kappa_m-1,\kappa_{m+1},\dots,\kappa_l)$. Then
\[
\DHom_{\hhh n}(\spek\la\kappa,\spek\mu\kappa)\cong\DHom_{\hhh{n-k}}(\spek{\la\bt}{\kappa\bt},\spek{\mu\bt}{\kappa\bt})
\]
as graded vector spaces over $\bbf$.
\end{cory}

\begin{pf}
{\abovedisplayskip=0pt\abovedisplayshortskip=0pt~\vspace*{-\baselineskip}
\begin{align*}
\DHom_{\hhh n}(\spek\la\kappa,\spek\mu\kappa)&\cong\DHom_{\hhh n}(\spek{\mu'}{\kappa'},\spek{\la'}{\kappa'})\tag*{ by \cref{homconj},}\\
&\cong\DHom_{\hhh{n-k}}(\spek{(\mu\bt)'}{(\kappa\bt)'},\spek{(\la\bt)'}{(\kappa\bt)'})\tag*{ by \cref{cr},}\\
&\cong\DHom_{\hhh{n-k}}(\spek{\la\bt}{\kappa\bt},\spek{\mu\bt}{\kappa\bt})\tag*{ by \cref{homconj} again.\qedhere}
\end{align*}}
\end{pf}

Now we prove a `final column removal' theorem, where we assume that the rightmost non-empty columns of $\la$ and $\mu$ are in the same place and of the same length.

\begin{thm}[Final Column Removal]\label{fcr}
Suppose $\la,\mu\in\mptn ln$ and $1\ls m\ls l$. Suppose $\la^{(1)}=\dots=\la^{(m-1)}=\mu^{(1)}=\dots=\mu^{(m-1)}=\varnothing$, $d:=\la^{(m)}_1=\mu^{(m)}_1$ and $k:=(\la^{(m)})'_d=(\mu^{(m)})'_d$. Let $\la\lf=\la\lf(d-1,m)$, $\mu\lf=\mu\lf(d-1,m)$ and $\kappa\lf=(\kappa_m,\dots,\kappa_l)$. Then
\[
\DHom_{\hhh n}(\spek\la\kappa,\spek\mu\kappa)\cong\DHom_{\hhh{n-k}}(\spek{\la\lf}{\kappa\lf},\spek{\mu\lf}{\kappa\lf})
\]
as graded vector spaces over $\bbf$.
\end{thm}

\begin{pf}
We first use \cref{rr} to remove the first $k$ rows of length $d$ from both $\la^{(m)}$ and $\mu^{(m)}$. We obtain
\begin{align*}
\DHom_{\hhh n}(\spek\la\kappa,\spek\mu\kappa)&\cong\DHom_{\hhh{n-dk}}(\spek{\la\bt}{\kappa\bt},\spek{\mu\bt}{\kappa\bt})
\\
\intertext{where $\la\bt=\la\bt(k,m)$, $\mu\bt=\mu\bt(k,m)$ and $\kappa\bt=(\kappa_m-k,\kappa_2,\dots,\kappa_l)$. We then use \cref{rr} again to add $k$ rows of length $d-1$ to the top of both $\la\bt^{(m)}$ and $\mu\bt^{(m)}$. We obtain}
\DHom_{\hhh{n-dk}}(\spek{\la\bt}{\kappa\bt},\spek{\mu\bt}{\kappa\bt})&\cong\DHom_{\hhh{n-k}}(\spek{\la\lf}{\kappa\lf},\spek{\mu\lf}{\kappa\lf})
\end{align*}
which gives the result.
\end{pf}

It will be helpful below to be able to give a direct construction for final column removal, as done in the proof of \cref{cr} for first column removal. We assume the hypotheses and notation of \cref{fcr}, and for ease of notation we assume that $\spe\la$ and $\spe\mu$ are defined using the \emc $\kappa$, while $\spe{\la\lf}$ and $\spe{\mu\lf}$ are defined using $\kappa\lf$. We can also assume that $\cont(\la)=\cont(\mu)=:\alpha$, and hence $\cont(\la\lf)=\cont(\mu\lf)=:\beta$.

We identify $\sss{n-k}$ with its image under the map $\shift0:\sss{n-k}\to\sss n$, and similarly for $\hhh\beta$ and $\hhh\alpha$. As in the proof of \cref{cr} we make an assumption on preferred reduced expressions: given a standard $\mu\lf$-tableau $\ttt$, we define $\ttt^+$ to be the standard $\mu$-tableau obtained by adding a column with entries $n-k+1,\dots,n$ at the right of component $m$; then we have $w_{\ttt^+}=w_\ttt$, and we assume that our preferred reduced expressions have been chosen in such a way that $\psi_{\ttt^+}=\psi_\ttt$.

\begin{lemma}\label{annint}
With the above notation, we have $\ann(z_{\la\lf})=\ann(z_\la)\cap\hhh\beta$.
\end{lemma}

\begin{pf}
It follows directly from the presentation for column Specht modules that $\ann(z_{\la\lf})\subseteq\ann(z_\la)\cap\hhh\beta$, so we must show the opposite containment. Consider the $\hhh\beta$-submodule $\hhh\beta z_\la$ of $\spe\la$ generated by $z_\la$. For any $\ttt\in\std{\la\lf}$ we have $v_{\ttt^+}=\psi_{\ttt^+}z_\la=\psi_\ttt z_\la\in\hhh\beta z_\la$, and the $v_{\ttt^+}$ are linearly independent, so $\dim_\bbf\hhh\beta z_\la\gs|\std{\la\lf}|=\dim_\bbf\spe{\la\lf}$. So we have
\begin{align*}
\dim_\bbf\hhh\beta z_\la&\gs\dim_\bbf\hhh\beta z_{\la\lf},\\
\intertext{i.e.}
\dim_\bbf\frac{\hhh\beta}{\ann(z_\la)\cap\hhh\beta}&\gs\dim_\bbf\frac{\hhh\beta}{\ann(z_{\la\lf})},
\end{align*}
and so $\ann(z_{\la\lf})\supseteq\ann(z_\la)\cap\hhh\beta$.
\end{pf}

Now we consider dominated homomorphisms. Observe that since $\la$ and $\mu$ have the same last column, $\cstd\la\mu=\rset{\ttt^+}{\ttt\in\cstd{\la\lf}{\mu\lf}}$. So if $\phi\in\chom\la\mu$, then we can write
\[
\phi(z_\la)=\sum_{\ttt\in\cstd{\la\lf}{\mu\lf}}a_\ttt v_{\ttt^+}\quad\text{with }a_\ttt\in\bbf.
\]
Then we can define a homomorphism
\begin{align*}
\phi^-:\spek{\la\lf}{\kappa\lf}&\longrightarrow\spek{\mu\lf}{\kappa\lf}\\
z_{\la\lf}&\longmapsto\sum_{\ttt\in\cstd{\la\lf}{\mu\lf}}a_\ttt v_\ttt.
\end{align*}
To see that this definition yields a well-defined homomorphism, we must show that $h\sum_\ttt a_\ttt v_\ttt=0$ whenever $h\in\ann(z_{\la\lf})$. By \cref{annint} we have $h\in\ann(z_\la)$, and hence (since $\phi$ is a homo\-morphism) $h\sum_\ttt a_\ttt v_{\ttt^+}=0$; in other words, $h\sum_{\ttt}a_\ttt\psi_\ttt\in\ann(z_\mu)$. We also have $h\sum_\ttt a_\ttt\psi_\ttt\in\hhh\beta$, so by \cref{annint} again (with $\la$ replaced by $\mu$) $h\sum_\ttt a_\ttt\psi_\ttt\in\ann(z_{\mu\lf})$, as required.

So we have a map $\phi\mapsto\phi^-:\dhom_{\hhh n}(\spek\la\kappa,\spek\mu\kappa)\to\dhom_{\hhh{n-k}}(\spek{\la\lf}{\kappa\lf},\spek{\mu\lf}{\kappa\lf})$. This is obviously an injective map of degree 0, and hence (by \cref{fcr}) a graded isomorphism.

\subsection{Generalised column removal}\label{gcrsec}

Armed with first column removal and final column removal, we can now consider generalised column removal. In what follows, we fix $c\gs0$ and $1\ls m\ls l$, and for any $\nu\in\mptn ln$ we write $\nu\lf=\nu\lf(c,m)$ and $\nu\rg=\nu\rg(c,m)$. We suppose $\la,\mu\in\mptn ln$, and assume that $|\la\lf|=|\mu\lf|=:n\lf$, so that $|\la\rg|=|\mu\rg|=n-n\lf=:n\rg$. We also assume that $\la\dom\mu$. This assumption implies that $\la\lf\dom\mu\lf$ and $\la\rg\dom\mu\rg$, which in particular gives
\[
(\la^{(m)})'_c\gs(\mu^{(m)})'_c\gs(\mu^{(m)})'_{c+1}
\]
so that it is possible to define a multipartition $\lr\la\mu\in\mptn ln$ with $(\lr\la\mu)\lf=\la\lf$ and $(\lr\la\mu)\rg=\mu\rg$.

We write $\kappa\lf=(\kappa_m,\dots,\kappa_l)$, $\kappa\rg=(\kappa_1,\dots,\kappa_m+c)$, $\scrh\lf=\hhh{n\lf}$ and $\scrh\rg=\hhh{n\rg}$. For ease of notation, we will assume throughout the following that the Specht modules $\spe\la$, $\spe\mu$ and $\spe{\lr\la\mu}$ are defined using the \emc $\kappa$, while $\spe{\la\lf}$ and $\spe{\mu\lf}$ are defined using $\kappa\lf$ and $\spe{\la\rg}$ and $\spe{\mu\rg}$ are defined using $\kappa\rg$.

Suppose $\phi\lf\in\dhom_{\scrh\lf}(\spe{\la\lf},\spe{\mu\lf})$ and $\phi\rg\in\dhom_{\scrh\rg}(\spe{\la\rg},\spe{\mu\rg})$, and write
\[
\phi\lf(z_{\la\lf})=\sum_{\tts\in\cstd{\la\lf}{\mu\lf}}a_\tts v_\tts,\qquad
\phi\rg(z_{\la\rg})=\sum_{\ttt\in\cstd{\la\rg}{\mu\rg}}b_\ttt v_\ttt
\]
with coefficients $a_\tts,b_\ttt\in\bbf$. If there is a homomorphism $\phi:\spe\la\to\spe\mu$ satisfying
\[
\phi(z_\la)=\sum_{\substack{\tts\in\cstd{\la\lf}{\mu\lf}\\\ttt\in\cstd{\la\rg}{\mu\rg}}}a_\tts b_\ttt v_{\tts\cj\ttt},
\]
then we write $\phi=\phi\lf\cj\phi\rg$, and say that $\phi$ is a \emph{product homomorphism}.

\begin{lemma}\label{factor}
Every product homomorphism $\spe\la\to\spe\mu$ factors through $\spe{\lr\la\mu}$.
\end{lemma}

\begin{pf}
Suppose that $\phi=\lr\phi\phi$ is a product homomorphism, and as above write
\[
\phi\lf(z_{\la\lf})=\sum_{\tts\in\cstd{\la\lf}{\mu\lf}}a_\tts v_\tts,\qquad
\phi\rg(z_{\la\rg})=\sum_{\ttt\in\cstd{\la\rg}{\mu\rg}}b_\ttt v_\ttt.
\]
Now define
\begin{alignat*}2
\phi\lf\cj\id:\spe{\lr\la\mu}&\longrightarrow\spe\mu&\id\cj\phi\rg:\spe\la&\longrightarrow\spe{\lr\la\mu}\\
z_{\lr\la\mu}&\longmapsto\sum_{\tts\in\cstd{\la\lf}{\mu\lf}}a_\tts v_{\tts\cj\ttt_{\mu\rg}}\qquad&z_\la&\longmapsto\sum_{\ttt\in\cstd{\la\rg}{\mu\rg}}b_\ttt v_{\ttt_{\la\lf}\cj\ttt}.
\end{alignat*}
Then $\phi\lf\cj\id$ and $\id\cj\phi\rg$ are both $\hhh n$-homomorphisms; this follows from the direct constructions of column removal homomorphisms in the proof of \cref{cr} and following the proof of \cref{fcr}. Clearly $(\phi\lf\cj\id)\circ(\id\cj\phi\rg)=\phi$, so $\phi$ factors through $\spe{\lr\la\mu}$.
\end{pf}

\begin{propn}\label{sumpro}
Assume the hypotheses and notation above. Then every $\phi\in\DHom_{\hhh n}(\spe\la,\spe\mu)$ is a sum of product homomorphisms.
\end{propn}

\begin{pf}
We may assume that $\cont(\la)=\cont(\mu)$ (since otherwise there are no non-zero homomorphisms $\spe\la\to\spe\mu$). So for this proof we write $\alpha:=\cont(\la)$ and define $\shift0$ to be the map from $\scrh\lf$ to $\hhh\alpha$ obtained by combining the maps $\shift0:\hhh\beta\to\hhh\alpha$ for all positive roots $\beta$ of height $n\lf$; similarly, $\shift{n\lf}$ denotes the map from $\scrh\rg$ to $\hhh\alpha$ obtained by combining the maps $\shift{n\lf}:\hhh\beta\to\hhh\alpha$ for all $\beta$ of height $n\rg$.

For this proof we make an assumption about the choice of preferred reduced expressions similar to that in the proof of \cref{fcr}. Specifically, we assume that these expressions have been chosen in such a way that if $\tts\in\cstd{\la\lf}{\mu\lf}$ and $\ttt\in\cstd{\la\rg}{\mu\rg}$, then the preferred expression for $w_{\tts\cj\ttt}$ is just the concatenation of the preferred expression for $w_\tts$ with the expression obtained by applying $\shift{n\lf}$ to every term in the preferred expression for $w_\ttt$. Hence $\psi_{\tts\cj\ttt}=\psi_\tts\shift{n\lf}(\psi_\ttt)$.

Now we show that every dominated homomorphism $\spe\la\to\spe\mu$ is a sum of product homomorphisms. To do this, we first discuss dominated tableaux. Note that the conditions on $\la$ and $\mu$ imply that $\cstd\la\mu=\rset{\lr\ttt\ttt}{\ttt\lf\in\cstd{\la\lf}{\mu\lf},\ \ttt\rg\in\cstd{\la\rg}{\mu\rg}}.$ Choose a total order $\blacktriangleright$ on $\cstd\la\mu$ with the property that if $\ttr,\tts\in\cstd{\la\lf}{\mu\lf}$ and $\ttt,\ttu\in\cstd{\la\rg}{\mu\rg}$, then
\[
\ttr\cj\ttt\blacktriangleright\ttr\cj\ttu\Longleftrightarrow\tts\cj\ttt\blacktriangleright\tts\cj\ttu\qquad\text{and}\qquad
\ttr\cj\ttt\blacktriangleright\tts\cj\ttt\Longleftrightarrow\ttr\cj\ttu\blacktriangleright\tts\cj\ttu.
\]
(For example, we could do this by choosing total orders $\blacktriangleright\lf,\blacktriangleright\rg$ on $\cstd{\la\lf}{\mu\lf},\cstd{\la\rg}{\mu\rg}$ and setting $\ttv\blacktriangleright\ttw$ if and only if $\ttv\lf\blacktriangleright\lf\ttw\lf$ or ($\ttv\lf=\ttw\lf$ and $\ttv\rg\blacktriangleright\rg\ttw\rg$).)

Now suppose $\varphi:\spe\la\to\spe\mu$ is a non-zero dominated homomorphism, and write $\varphi(z_\la)=\sum_{\ttt\in\cstd\la\mu} a_\ttt v_\ttt$ with each $a_\ttt\in\bbf$. Let $\ttu$ be the largest tableau (with respect to $\blacktriangleright$) such that $a_\ttu\neq0$, and proceed by induction on $\ttu$.

\clam
Let $\calu$ denote the set of tableaux $\ttt\in\cstd\la\mu$ such that $\ttt\rg=\ttu\rg$. Then there is an $\scrh\lf$-homomorphism
\begin{align*}
\varphi\lf^\ttu:\spe{\la\lf}&\longrightarrow\spe{\mu\lf}\\
z_{\la\lf}&\longmapsto\sum_{T\in\calu}a_\ttt v_{\ttt\lf}.
\end{align*}
\prof
First we make an observation, which follows from the construction of Specht modules and our assumptions on preferred reduced expressions. If $\ttw\in\cstd\la\mu$ and $h\in\scrh\lf$, and we write $hv_{\ttw\lf}=\sum_{\ttt\in\std{\mu\lf}}b_\ttt v_\ttt$, then $\shift0(h)v_\ttw=\sum_{\ttt\in\std{\mu\lf}}b_\ttt v_{\ttt\cj\ttw\rg}$. In particular, $\shift0(h)v_\ttw$ is a linear combination of basis elements $v_\tts$ for $\tts\in\stdlr\mu$ with $\tts\rg=\ttw\rg$.

Now take $h\in\ann(z_{\la\lf})$. Then $\shift0(h)\in\ann(z_\la)$, so $\shift0(h)\sum_{\ttt\in\cstd\la\mu}a_\ttt v_\ttt=0$ (because $\varphi$ is a homomorphism). If we look just at $\shift0(h)\sum_{\ttt\in\calu}a_\ttt v_\ttt$, then by the previous paragraph this lies in $\lspan{v_\ttt}{\ttt\in\stdlr\mu,\ \ttt\rg=\ttu\rg}_\bbf$, while $\shift0(h)\sum_{\ttt\notin\calu}a_\ttt v_\ttt$ lies in $\lspan{v_\ttt}{\ttt\in\stdlr\mu,\ \ttt\rg\neq\ttu\rg}_\bbf$. The $v_\ttt$ are linearly independent, and hence
\[
\lspan{v_\ttt}{\ttt\in\stdlr\mu,\ \ttt\rg=\ttu\rg}_\bbf\cap\lspan{v_\ttt}{\ttt\in\stdlr\mu,\ \ttt\rg\neq\ttu\rg}_\bbf=0.
\]
Hence $\shift0(h)\sum_{\ttt\in\calu}a_\ttt v_\ttt=0$.

Define a linear map $\cj\ttu\rg:\spe{\mu\lf}\to\spe\mu$ by $v_\ttt\mapsto v_{\ttt\cj\ttu\rg}$ for $\ttt\in\std{\mu\lf}$. Then, from above, we have
\[
(hm)\cj\ttu\rg = h(m\cj\ttu\rg)
\]
for any $h\in\scrh\lf$ and any $m\in\spe{\mu\lf}$. So for each $h\in\ann(z_{\la\lf})$, we have $h\sum_{\ttt\in\calu}a_\ttt v_{\ttt\lf}=0$.
\malc
We can do essentially the same thing left and right interchanged; that is, if we let $\calu'=\lset{\ttt\in\cstd\la\mu}{\ttt\lf=\ttu\lf}$, then we have an $\scrh\rg$-homomorphism
\begin{align*}
\varphi\rg^\ttu:\spe{\la\rg}&\longrightarrow\spe{\mu\rg}\\
z_{\la\rg}&\longmapsto\sum_{\ttt\in\calu'}a_\ttt v_{\ttt\rg}.
\end{align*}

As in the proof of \cref{factor}, we construct homomorphisms
\[
\phi\lf^\ttu\cj\id:\spe{\lr\la\mu}\longrightarrow\spe\mu\qquad\id\cj\phi\rg^\ttu:\spe\la\longrightarrow\spe{\lr\la\mu},
\]
whose composition is the product homomorphism $\lr{\phi^\ttu}{\phi^\ttu}:\spe\la\to\spe\mu$. $v_\ttu$ appears with non-zero coefficient (namely $a_\ttu^2$) in $\lr{\phi^\ttu}{\phi^\ttu}$, and $\ttu$ is maximal (with respect to the order $\blacktriangleright$) with this property. So if we consider the homomorphism $\xi:=\varphi-\mfrac1{a_\ttu}\lr\phi\phi$, then (if $\xi\neq0$) the most dominant tableau occurring with non-zero coefficient in $\xi$ is smaller than $\ttu$. By induction $\xi$ is a sum of product homomorphisms, and hence so is $\varphi$.
\end{pf}

Now we can prove our main result.

\begin{thm}[Generalised graded column removal]\label{gcr}
Suppose $\la,\mu\in\mptn ln$, $c\gs0$ and $1\ls m\ls l$ and define $\la\lf,\la\rg,\mu\lf,\mu\rg$ as in \cref{gcrnotation}. Assume $|\la\lf(c,m)|=|\mu\lf(c,m)|=:n\lf$ and $|\la\rg(c,m)|=|\mu\rg(c,m)|=:n\rg$ for some fixed $c\geq0$ and $1\leq m\leq l$ and define $\scrh\lf=\hhh{n\lf}$ and $\scrh\rg=\hhh{n\rg}$.
\begin{enumerate}
\item
For any $\phi\lf\in\dhom_{\scrh\lf}(\spe{\la\lf},\spe{\mu\lf})$ and $\phi\rg\in\dhom_{\scrh\rg}(\spe{\la\rg},\spe{\mu\rg})$, there is a product homomorphism $\lr\phi\phi\in\chom\la\mu$.
\item
The map $\phi\lf\otimes\phi\rg\mapsto\lr\phi\phi$ defines an isomorphism of graded $\bbf$-vector spaces
\[
\dhom_{\scrh\lf}(\spe{\la\lf},\spe{\mu\lf})\otimes\dhom_{\scrh\rg}(\spe{\la\rg},\spe{\mu\rg})\cong\chom\la\mu.
\]
\end{enumerate}
\end{thm}

\begin{pf}
First suppose $\la\ndom\mu$. Then $\cstd\la\mu=\emptyset$, so $\chom\la\mu=0$. Furthermore, we have either $\la\lf\ndom\mu\lf$ or $\la\rg\ndom\mu\rg$, so that either $\dhom_{\scrh\lf}(\spe{\la\lf},\spe{\mu\lf})=0$ or $\dhom_{\scrh\rg}(\spe{\la\rg},\spe{\mu\rg})=0$. So the result follows.

So we can assume that $\la\dom\mu$, which allows us to define the multipartition $\lr\la\mu$ as above. Applying \cref{cr} repeatedly, we have
\[
\DHom_{\hhh n}(\spe\la,\spe{\lr\la\mu})\cong\DHom_{\scrh\rg}(\spe{\la\rg},\spe{\mu\rg}).
\]
Similarly, by \cref{fcr} applied repeatedly we have
\[
\DHom_{\hhh n}(\spe{\lr\la\mu},\spe\mu)\cong\DHom_{\scrh\lf}(\spe{\la\lf},\spe{\mu\lf}).
\]
Combining these isomorphisms, and using the explicit constructions given above, we have an isomorphism of graded vector spaces
\begin{align*}
\DHom_{\scrh\lf}(\spe{\la\lf},\spe{\mu\lf})\otimes\DHom_{\scrh\rg}(\spe{\la\rg},\spe{\mu\rg})&\stackrel\sim\longrightarrow\DHom_{\hhh n}(\spe{\lr\la\mu},\spe\mu)\otimes\DHom_{\hhh n}(\spe\la,\spe{\lr\la\mu})\\
\phi\lf\otimes\phi\rg&\longmapsto(\phi\lf\cj\id)\otimes(\id\cj\phi\rg).
\end{align*}
Composition of homomorphisms yields a map
\[
\omega:\DHom_{\hhh n}(\spe{\lr\la\mu},\spe\mu)\otimes\DHom_{\hhh n}(\spe\la,\spe{\lr\la\mu})\longrightarrow\chom\la\mu
\]
which is homogeneous of degree zero, and by \cref{factor,sumpro} $\omega$ is surjective. So we have a surjective map
\begin{align*}
\DHom_{\scrh\lf}(\spe{\la\lf},\spe{\mu\lf})\otimes\DHom_{\scrh\rg}(\spe{\la\rg},\spe{\mu\rg})&\longrightarrow\chom\la\mu\\
\phi\lf\otimes\phi\rg&\longmapsto\lr\phi\phi.
\end{align*}
This map is easily seen to be injective, and the result follows. 
\end{pf}

\subsection{Generalised row removal}

Now we consider generalised row removal for homomorphisms between column Specht modules. Fix $1\ls m\ls l$ and $r\gs0$, and for any $\nu\in\mptn ln$ write $\nu\tp=\nu\tp(r,m)$, $\nu\bt=\nu\bt(r,m)$. Suppose $\la,\mu\in\mptn ln$ with $|\la\tp|=|\mu\tp|=:n\tp$, so that $|\la\bt|=|\mu\bt|=n-n\tp=:n\bt$. Set $\kappa\tp=(\kappa_1,\dots,\kappa_m)$ and $\kappa\bt=(\kappa_m-r,\kappa_{m+1},\dots,\kappa_l)$, and write $\scrh\tp=\hhh{n\tp}$ and $\scrh\bt=\hhh{n\bt}$. In what follows we shall take $\spe\la$ and $\spe\mu$ to be defined with respect to the \emc $\kappa$, $\spe{\la\tp}$ and $\spe{\mu\tp}$ with respect to $\kappa\tp$, and $\spe{\la\bt}$ and $\spe{\mu\bt}$ with respect to $\kappa\bt$.

With this notation in place, we can state a generalised row-removal theorem for homomorphisms.  This follows from \cref{gcr} using \cref{homconj} in the same way that \cref{rr} is deduced from \cref{cr}.

\begin{thm}[Generalised graded row removal]\label{grr}
Suppose $\la,\mu\in\mptn ln$, $r\gs0$ and $1\ls m\ls l$ and define $\la\tp,\la\bt,\mu\tp,\mu\bt,n\tp,n\bt,\scrh\tp,\scrh\bt$ as above. Then there is an isomorphism of graded $\bbf$-vector spaces
\[
\dhom_{\scrh\tp}(\spe{\la\tp},\spe{\mu\tp})\otimes\dhom_{\scrh\bt}(\spe{\la\bt},\spe{\mu\bt})\cong\chom\la\mu.
\]
\end{thm}

Our proof of \cref{gcr} gives a direct construction of the column-removal isomorphism, but a direct construction for row removal seems to be hard to obtain, especially using the standard bases for column Specht modules.

\begin{eg}
Take $e=2$ and $\kappa=(0,1,0)$. Let $\la=\big((1^2),(2,1^3),(1)\big)$ and $\mu=\big((1),(3,1),(3)\big)$, and take $(m,r)=(2,1)$, so that $\la\tp=\big((1^2),(2)\big)$, $\la\bt=\big((1^3),(1)\big)$ and $\mu\tp=\mu\bt=\big((1),(3)\big)$. Set $\kappa\tp=(0,1)$ and $\kappa\bt=(0,0)$. Then (regardless of the field $\bbf$) the graded dimensions of $\dhom_{\hhh4}(\spek{\la\tp}{\kappa\tp},\spek{\mu\tp}{\kappa\tp})$ and $\dhom_{\hhh4}(\spek{\la\bt}{\kappa\bt},\spek{\mu\bt}{\kappa\bt})$ are $v$ and $1$ respectively. So by \cref{grr} the graded dimension of $\dhom_{\hhh8}(\spek\la\kappa,\spek\mu\kappa)$ is $v$. The unique (up to scaling) homomorphisms
\begin{alignat*}3
\spe{\la\tp}&\longrightarrow\spe{\mu\tp},&\qquad\spe{\la\bt}&\longrightarrow\spe{\mu\bt},&\qquad\spe\la&\longrightarrow\spe\mu\\
\intertext{are given by}
z_{\la\tp}&\longmapsto v_\tts,&z_{\la\bt}&\longmapsto v_\ttt,&z_\la&\longmapsto v_\ttu+2v_\ttv,
\end{alignat*}
where
\[
\tts=\gyoung(^3^\hf;3,/\hf,;1;2;4),\qquad\ttt=\gyoung(^3^\hf;2,/\hf,;1;3;4),\qquad\ttu=\gyoung(^7;7,/\hf,^3^\hf;2;6;8,^3^\hf;3,/\hf,;1;4;5),\quad\ttv=\gyoung(^7;7,/\hf,^3^\hf;4;6;8,^3^\hf;5,/\hf,;1;2;3).
\]
It seems hard to reconcile these homomorphisms when expressed in this form, except perhaps in characteristic $2$. (Note that the incompatibility of these expressions is not an artefact of the choice of preferred reduced expressions -- the standard basis elements appearing in this example are independent of the choice of reduced expressions.)
\end{eg}

In order to obtain an explicit row-removal construction, it seems to be necessary to use a different basis for the Specht module. Suppose we have $\la\bt$, $\la\tp$, $\mu\bt$ and $\mu\tp$ as above, with $|\mu\tp|=n\tp=|\la\tp|$. Partition the set $\{1,\dots,n\}$ into two sets $S\bt$ and $S\tp$, by defining $S\bt$ to be the set of integers in the bottom part of $\ttt_\la$ and $S\tp$ the set of integers in the top part; that is,
\begin{align*}
S\bt&=\lset{\ttt_\la(s,c,k)}{(s,c,k)\in[\la]\text{ and either }k>m\text{ or }k=m\text{ and }s>r},\\
S\tp&=\lset{\ttt_\la(s,c,k)}{(s,c,k)\in[\la]\text{ and either }k<m\text{ or }k=m\text{ and }s\ls r}.
\end{align*}
Let $\lab\bt:\{1,\dots,n\bt\}\rightarrow S\bt$ and $\lab\tp:\{1,\dots,n\tp\}\rightarrow S\tp$ be the unique order-preserving bijections.

Now given a $\mu\bt$-tableau $\ttt$ and a $\mu\tp$-tableau $\tts$, define a $\mu$-tableau $\ttt\rj\tts$ by composing $\lab\bt$ with $\ttt$ and $\lab\tp$ with $\tts$ and `gluing' in the natural way.

\begin{lemma}\label{joinstd}
Suppose $\la$ and $\mu$ satisfy the conditions above. If $\ttt\in\cstd{\la\bt}{\mu\bt}$ and $\tts\in\cstd{\la\tp}{\mu\tp}$, then $\ttt\rj\tts\in\cstd\la\mu$.
\end{lemma}

\begin{pf}
First we show that $\ttt\rj\tts$ is standard. Suppose $A$ and $B$ are nodes in the same component of $[\mu]$, with $B$ either immediately to the right of $A$ or immediately below $A$; then we require $\ttt\rj\tts(B)>\ttt\rj\tts(A)$. This is clear from the fact that $\tts$ and $\ttt$ are standard and the functions $\lab\tp$ and $\lab\bt$ are order-preserving, except in the case where $A=(r,b,m)$ and $B=(r+1,b,m)$ for some $1\ls b\ls\mu^{(m)}_{r+1}$. So assume we are in this situation.

Let $k=\la^{(m)}_{r+1}$. Then the first $k$ columns of $\la\tp^{(m)}$ all have length $r$. Since $\cstd{\la\tp}{\mu\tp}$ is non-empty we have $\la\tp\dom\mu\tp$, and hence the first $k$ columns of $\mu\tp^{(m)}$ all have length $r$ also. Hence (since $\tts$ is $\la\tp$-dominated) $\tts$ agrees with $\ttt_{\la\tp}$ on these columns. So we have $\ttt\rj\tts(A)=\lab\tp(\ttt_{\la\tp}(A))=\ttt_\la(A)$.

We also have $\la\bt\dom\mu\bt$ since $\cstd{\la\bt}{\mu\bt}\neq\emptyset$, so that $k\gs\mu^{(m)}_{r+1}\gs b$ (and in particular $B\in[\la]$). Since $\ttt$ is $\la\bt$-dominated, we have $\ttt(1,b,1)\gs\ttt_{\la\bt}(1,b,1)$, so that
\[
\ttt\rj\tts(B)=\lab\bt(\ttt(1,b,1))\gs\lab\bt(\ttt_{\la\bt}(1,b,1))=\ttt_\la(B).
\]
So $\ttt\rj\tts(A)=\ttt_\la(A)<\ttt_\la(B)\ls\ttt\rj\tts(B)$, as required.

To see that $\ttt\rj\tts$ is $\la$-dominated, it suffices to note that since $\tts\in\cstd{\la\tp}{\mu\tp}$, every element of $S\tp$ appears in $\lab\tp(\tts)$ at least as far to the left as it appears in $\ttt_\la$, and likewise for $\ttt\in\cstd{\la\bt}{\mu\bt}$ and elements of $S\bt$.
\end{pf}

Now we can give a conjectured explicit construction for the generalised row-removal isomorphism for homomorphisms. Recall from \cref{dualsec} the basis $\rset{f_\ttt}{\ttt\in\std\mu}$ for $(\rspe\mu)\dual$; using \cref{7.25} and shifting the degree of each $f_\ttt$ by $\df(\mu)$, we can regard $\rset{f_\ttt}{\ttt\in\std\mu}$ as a basis for $\spe\mu$. Note that by the analogue of \cref{triang}(\ref{triang2}) for column Specht modules, any $\phi\in\chom\la\mu$ can be written as
\[
\phi(z_\la)=\sum_{\ttt\in\cstd\la\mu} a_\ttt f_\ttt \quad \text{ for some } a_\ttt\in\bbf.
\]

\begin{conj}\label{exprow}
Suppose $\la,\mu\in\mptn ln$, $r\gs0$ and $1\ls m\ls n$. Define $\la\tp,\la\bt,\mu\tp,\mu\bt,n\tp,n\bt,\scrh\tp,\scrh\bt$ as above, and assume $|\mu\tp|=n\tp$. Suppose $\phi\tp\in\dhom_{\scrh\tp}(\spe{\la\tp},\spe{\mu\tp})$ and $\phi\bt\in\dhom_{\scrh\bt}(\spe{\la\bt},\spe{\mu\bt})$, and write
\[
\phi\bt(z_{\la\bt})=\sum_{\mathclap{\ttt\in\cstd{\la\bt}{\mu\bt}}}a_\ttt f_\ttt,\qquad \phi\tp(z_{\la\tp})=\sum_{\mathclap{\tts\in\cstd{\la\tp}{\mu\tp}}}b_\tts f_\tts
\]
with $a_\ttt,b_\tts\in\bbf$. Then there is an $\hhh n$-homomorphism $\phi\bt\rj\phi\tp:\spe\la\to\spe\mu$ satisfying
\[
\phi\bt\rj\phi\tp(z_\la)=\sum_{\mathclap{\substack{\ttt\in\cstd{\la\bt}{\mu\bt}\\ \tts\in\cstd{\la\tp}{\mu\tp}}}} a_\ttt b_\tts f_{\ttt\rj\tts}.
\]
\end{conj}

\begin{eg}
Retaining the notation from the last example, we have
\[
\ttt_\la=\gyoung(^4;7,^4;8,/\hf,:^\hf;2;6,:^\hf;3,:^\hf;4,:^\hf;5,/\hf,;1),
\]
so that $S\tp=\{2,6,7,8\}$ and $S\bt=\{1,3,4,5\}$. Taking $\tts$, $\ttt$ and $\ttu$ as in the last example, we get $\ttt\rj\tts=\ttu$. It is easy to check that
\[
f_\tts=v_\tts,\qquad f_\ttt=v_\ttt,\qquad f_\ttu=v_\ttu+2v_\ttv,
\]
so the conjecture holds in this case.
\end{eg}

\begin{rmk}
If \cref{exprow} is true, then we have a map of graded $\bbf$-vector spaces
\begin{align*}
\dhom_{\scrh\bt}(\spe{\la\bt},\spe{\mu\bt})\otimes\dhom_{\scrh\tp}(\spe{\la\tp},\spe{\mu\tp})&\longrightarrow\chom\la\mu\\
\phi\bt\otimes\phi\tp&\longmapsto\phi\bt\rj\phi\tp.
\end{align*}
This map is obviously linear, and (since the $f_\ttt$ are linearly independent) injective. Hence by \cref{grr} it is a bijection. So we have an explicit construction for the generalised row-removal isomorphism.
\end{rmk}

\section{Index of notation}\label{indsec}

For the reader's convenience we conclude with an index of the notation we use in this paper. We provide references to the relevant subsections.

\newlength\colwi
\newlength\colwii
\newlength\colwiii
\setlength\colwi{2.5cm}
\setlength\colwiii{1cm}
\setlength\colwii\textwidth
\addtolength\colwii{-\colwi}
\addtolength\colwii{-\colwiii}
\addtolength\colwii{-1em}
\begin{longtable}{@{}p{\colwi}p{\colwii}p{\colwiii}@{}}
$\bbf$&a field&\\
$\bbn$&the set of positive integers&\\
$\sss n$&the symmetric group of degree $n$&\ref{symsec}\\
$s_1,\dots,s_{n-1}$&the Coxeter generators of $\sss n$&\ref{symsec}\\
$l$&the Coxeter length function on $\sss n$&\ref{symsec}\\
$\ls\Lf$&the left order on $\sss n$&\ref{symsec}\\
$\cls$&the Bruhat order on $\sss n$&\ref{symsec}\\
$\shift k$&the shift homomorphism $\sss m\to\sss n$&\ref{symsec}\\
$I$&the set $\zez$ (or $\bbz$, if $e=\infty$)&\ref{liesec}\\
$\Gamma$&a quiver with vertex set $I$&\ref{liesec}\\
$i\to j$&there is an arrow from $i$ to $j$ (but no arrow from $j$ to $i$) in $\Gamma$&\ref{liesec}\\
$i\rightleftarrows j$&there are arrows from $i$ to $j$ and from $j$ to $i$ in $\Gamma$&\ref{liesec}\\
$\alpha_i$&simple root labelled by $i\in I$&\ref{liesec}\\
$\La_i$&fundamental dominant weight labelled by $i\in I$&\ref{liesec}\\
$\bil\,\,$&invariant bilinear form&\ref{liesec}\\
$\La_\kappa$&the dominant weight $\La_{\kappa_1}+\dots+\La_{\kappa_l}$&\ref{liesec}\\
$\df(\alpha)$&$\bil{\La_k}\alpha-\frac12\bil\alpha\alpha$&\ref{liesec}\\
$\mptn ln$&the set of $l$-multipartitions of $n$&\ref{mptnsec}\\
$|\la|$&the number of nodes of a (multi)partition $\la$&\ref{mptnsec}\\
$\dom$&the dominance order on multipartitions or tableaux&\ref{mptnsec}\\
$[\la]$&the Young diagram of a multipartition $\la$&\ref{mptnsec}\\
$\varnothing$&the unique partition or $l$-multipartition of $0$&\ref{mptnsec}\\
$\la'$&the conjugate (multi)partition to $\la$&\ref{mptnsec}\\
$\std\la$&the set of standard $\la$-tableaux&\ref{mptnsec}\\
$\ttt'$&the conjugate tableau to $\ttt$&\ref{mptnsec}\\
$i\downarrow_\ttt j$&$i$ and $j$ lie in the same column of $\ttt$, with $j$ lower than $i$&\ref{mptnsec}\\
$i\swarrow_\ttt j$&$i$ and $j$ lie in the same component of $\ttt$, with $j$ strictly lower and to the left of $i$&\ref{mptnsec}\\
$i\Swarrow_\ttt j$&$i\swarrow_\ttt j$ or $i$ lies in an earlier component of $\ttt$ than $j$&\ref{mptnsec}\\
$\ttt_\la$&the $\la$-tableau obtained by writing $1,\dots,n$ in order down successive columns&\ref{mptnsec}\\
$\ttt^\la$&the $\la$-tableau obtained by writing $1,\dots,n$ in order along successive rows&\ref{mptnsec}\\
$w_\ttt$&the permutation for which $w_\ttt\ttt_\la=\ttt$&\ref{mptnsec}\\
$w^\ttt$&the permutation for which $w^\ttt\ttt^\la=\ttt$&\ref{mptnsec}\\
$\res A$&the residue of a node $A$&\ref{resdegsec}\\
$\cont(\la)$&the content of a multipartition $\la$&\ref{resdegsec}\\
$\df(\la)$&the defect of a multipartition $\la$&\ref{resdegsec}\\
$i(\ttt)$&the residue sequence of a tableau $\ttt$&\ref{resdegsec}\\
$i_\la$&$i(\ttt_\la)$&\ref{resdegsec}\\
$i^\la$&$i(\ttt^\la)$&\ref{resdegsec}\\
$\deg(\ttt)$&the degree of a tableau $\ttt$&\ref{resdegsec}\\
$\codeg(\ttt)$&the codegree of a tableau $\ttt$&\ref{resdegsec}\\
$\hhh n$&the KLR algebra of degree $n$&\ref{klrsec}\\
$\shift k$&the shift homomorphism $\hhh\beta\to\hhh\alpha$&\ref{klrsec}\\
$\hhh n^\kappa$&the cyclotomic KLR algebra determined by $\kappa$&\ref{klrsec}\\
$\bfB_A$&the Garnir belt corresponding to a Garnir node $A$&\ref{spechtsec}\\
$\garn A$&the Garnir element corresponding to a Garnir node $A$&\ref{spechtsec}\\
$\spe\la$&the column Specht module corresponding to a multipartition $\la$&\ref{spechtsec}\\
$\rspe\la$&the row Specht module corresponding to a multipartition $\la$&\ref{spechtsec}\\
$z_\la$&the standard generator of $\spe\la$&\ref{spechtsec}\\
$\psi_\ttt$&$\psi_{t_1}\dots\psi_{t_b}$, where $s_{t_1}\dots s_{t_b}$ is the preferred reduced expression for $w_\ttt$&\ref{spechtsec}\\
$v_\ttt$&$\psi_\ttt z_\la$&\ref{spechtsec}\\
$\shp\ttt m$&the $l$-multicomposition formed from the nodes of $\ttt$ whose entries are less than or equal to $m$&\ref{tabsec}\\
$\cstd\la\mu$&the set of $\la$-dominated standard $\mu$-tableaux&\ref{domtsec}\\
$\rstd\la\mu$&the set of $\la$-row-dominated standard $\mu$-tableaux&\ref{domtsec}\\
$\chom\la\mu$&the space of dominated homomorphisms from $\spe\la$ to $\spe\mu$&\ref{domhsec}\\
$\rhom\la\mu$&the space of dominated homomorphisms from $\rspe\la$ to $\rspe\mu$&\ref{domhsec}\\
$M\dual$&the graded dual of a graded module $M$&\ref{dualsec}\\
$M\lan k\ran$&the graded module $M$ with the grading shifted by $k$&\ref{dualsec}\\
$\stdlr\la$&the set of $\la$-tableaux in which the entries $1,\dots,n\lf$ appear strictly to the left of the entries $n\lf+1,\dots,n$&\ref{crhsec}\\
$\lr\la\la$&the multipartition obtained by joining the left and right parts $\la\lf,\la\rg$ together&\ref{gcrsec}\\
$\lr\ttt\ttt$&the tableau obtained by joining the left and right parts $\ttt\lf,\ttt\rg$ together&\ref{gcrsec}\\
\end{longtable}

\end{document}